\documentclass[aps,prd,amsmath,floats,floatfix,twocolumn,
  superscriptaddress,nofootinbib,showpacs]{revtex4-2}

\usepackage[dvipsnames]{xcolor}
\usepackage{hyperref}
\hypersetup{
  colorlinks,
  citecolor={cyan!50!black},
  linkcolor={orange!60!black},
  urlcolor=violet,
}
\usepackage{orcidlink}
\usepackage{graphicx}
\usepackage{amssymb}
\usepackage{mathtools}  %
\usepackage{tensor}
\usepackage{appendix}
\usepackage{bm} %
\usepackage{commath}  %
\usepackage{physics}
\usepackage{bbm} %
\usepackage{microtype}
\usepackage{enumitem}
\usepackage{siunitx}
\usepackage{csquotes}
\usepackage[capitalise]{cleveref}
\usepackage[caption=false]{subfig}
\newcommand{\ccite}[1]{Ref.~\cite{#1}}
\newcommand{\ccites}[1]{Refs.~\cite{#1}}
\newcommand{\Ccite}[1]{Reference~\cite{#1}}
\newcommand{\tikzinput}[1]{
  \includegraphics{#1.pdf}
}

\setlist[description]{font=\normalfont\itshape}

\usepackage[normalem]{ulem} %

\newcommand{\defeq}{\coloneqq}

\newcommand{\dgF}{\mathcal{F}}
\newcommand{\dgFi}[1]{\tensor{\dgF}{_{\!{#1}}^i}}
\newcommand{\dgFj}[1]{\tensor{\dgF}{_{\!{#1}}^j}}
\newcommand{\dgS}{\mathcal{S}}
\newcommand{\dgb}{b}
\newcommand{\dgf}{f}

\newcommand{\bas}{\psi}
\newcommand{\lagr}{\ell}
\newcommand{\jac}{\mathrm{J}}
\newcommand{\invjac}{{(\jac^{-1})}}
\newcommand{\surf}[1]{{#1}^\Sigma}

\newcommand{\grd}[1]{\underline{#1}}
\newcommand{\dgA}{\mathcal{A}}
\newcommand{\dgM}{M}
\newcommand{\dgD}{D}
\newcommand{\dgMD}{M\!\!D}
\newcommand{\dgL}{L}
\newcommand{\dgML}{M\!\!L}
\newcommand{\prol}{P}
\newcommand{\restr}{R}

\newcommand{\ident}{\mathbbm{1}}

\newcommand{\penaltyparam}{C}

\newcommand{\displ}{\xi}
\newcommand{\strain}{S}
\newcommand{\stress}{T}
\newcommand{\constrel}{Y}
\newcommand{\besselJ}{J}

\newcommand{\spec}{\texttt{SpEC}}
\newcommand{\spectre}{\texttt{SpECTRE}}

\begin{document}

\title{Unified discontinuous Galerkin scheme for a large class of elliptic equations}

\newcommand{\aei}{\affiliation{Max Planck Institute for Gravitational Physics
(Albert Einstein Institute), Am M{\"u}hlenberg 1, Potsdam 14476, Germany}}

\author{Nils L. Fischer\,\orcidlink{0000-0002-5767-3949}} \email{nils.fischer@aei.mpg.de} \aei
\author{Harald P. Pfeiffer\,\orcidlink{0000-0001-9288-519X}} \aei

\date{\today}

\begin{abstract}
  We present a discontinuous Galerkin internal-penalty scheme that is applicable
  to a large class of linear and nonlinear elliptic partial differential
  equations. The unified scheme can accommodate all second-order elliptic
  equations that can be formulated in first-order flux form, encompassing
  problems in linear elasticity, general relativity, and hydrodynamics, including
  problems formulated on a curved manifold. It allows for a wide range of linear
  and nonlinear boundary conditions, and accommodates curved and nonconforming
  meshes. Our generalized internal-penalty numerical flux and our
  Schur-complement strategy of eliminating auxiliary degrees of freedom make the
  scheme compact without requiring equation-specific modifications. We
  demonstrate the accuracy of the scheme for a suite of numerical test problems.
  The scheme is implemented in the open-source \spectre{} numerical relativity
  code.
\end{abstract}

\maketitle

\section{Introduction}\label{sec:intro}
    
Many problems in physics involve the numerical solution of second-order
elliptic partial differential equations (PDEs). Such elliptic problems often
represent static field configurations under the effect of external forces and
arise, for example, in electrodynamics, in linear or nonlinear elasticity, and
in general relativity. Elliptic problems also often accompany time evolutions,
where they constrain the evolved fields at every instant in time or provide
admissible initial data for the evolution.

Discontinuous Galerkin (DG) methods are gaining popularity in the computational
physics and engineering community and are currently most prevalently used for
time evolutions of hyperbolic boundary-value problems \cite{Reed1973-hw,
HesthavenWarburton, Cockburn2000-av, Cockburn2001-md}. Many properties that make
DG methods advantageous for time evolutions also apply
to elliptic problems, which lead to the development of DG schemes for elliptic
PDEs \cite{Wheeler1978-tj, Arnold2002}. In particular, DG schemes provide a
flexible mechanism for refining the computational grid, retaining exponential
convergence even in the presence of discontinuities when adaptive
mesh refinement (AMR) techniques are employed \cite{Schoetzau2014,
Vincent2019qpd}. Furthermore, some difficulties with DG schemes in time
evolutions, such as shock capturing, are not present in elliptic problems and
their static nature makes it often (but not always) straightforward to place
grid boundaries at discontinuities, thus relieving the AMR scheme from the
responsibility of resolving them. See, e.g., \ccite{HesthavenWarburton} and the
seminal paper~\cite{Arnold2002} for extensive discussions of DG schemes for the
Poisson equation, and \ccites{Hansbo2002-zb, Cockburn2006-rz, Ortner2007-wk} for
discussions of linear and nonlinear elasticity.

In the context of relativistic astrophysics and numerical relativity, DG methods
have been developed for hyperbolic equations on curved manifolds thus
far~\cite{Teukolsky2016-ja, Kidder2017-nz, Fambri2018-nu}. In
\ccite{Vincent2019qpd} we explored the feasibility of the DG method for elliptic
problems in numerical relativity confined to flat Poisson-type equations with
nonlinear sources. In this article we present a DG scheme suitable to solve a
significantly larger class of elliptic problems that arise in numerical
relativity. Most notably, the scheme encompasses the extended conformal thin
sandwich (XCTS) formulation of the general-relativistic Einstein constraint
equations on a curved manifold, and associated boundary conditions
\cite{Pfeiffer2004-oo, Cook2004-yf, BaumgarteShapiro}. Solutions to the XCTS
equations provide admissible initial data for general-relativistic time
evolutions, for scenarios such as two orbiting black holes or neutron stars
\cite{Lovelace2008-sw, Varma2018-fp, Tacik2015-iz, Tacik2016-me}. To
our knowledge, this article presents the first discontinuous Galerkin solution
of the full Einstein constraint equations. Aimed at applications in numerical
relativity, the scheme is implemented in the publicly available \spectre{} code
\cite{Kidder2017-nz, spectre, ellsolver}.

Furthermore, the elliptic DG scheme presented in this article is not limited to
applications in numerical relativity. It applies to all second-order elliptic
problems that can be formulated in first-order flux form. Besides the classic
Poisson and elasticity equations it covers a large class of elliptic problems in
general relativity and hydrodynamics, including coupled systems of equations and
those formulated on a curved manifold. With our unified DG scheme, new elliptic
systems can be implemented by supplying their first-order fluxes and sources,
hence no knowledge of the DG technology or of finite-element formulations is
required. This lowers the barrier for extending the capabilities of a simulation
code. We pay particular attention to support a wide range of linear and
nonlinear boundary conditions so our DG scheme is suited to solve many
real-world scenarios (as well as some out-of-this-world scenarios such as
initial data for evolutions of black holes and neutron stars). We are aware only of
\ccite{Feistauer2019-cw} studying a nonlinear boundary condition for an
elliptic DG problem.

To formulate the unified DG scheme we present a generalized internal-penalty
numerical flux, which avoids problem-specific parameters that are needed, e.g.,
in \ccites{Hartmann2008-zi, Epshteyn2007-dm, De_Basabe2008-wf}. We eliminate
auxiliary degrees of freedom that arise from the first-order form with a
Schur-complement strategy, which has proven more suitable to the unified DG
scheme than primal formulations that are commonly employed in the
literature~\cite{HesthavenWarburton, Arnold2002}. The resulting DG scheme is
compact, in the sense that it involves only nearest-neighbor couplings and no
auxiliary degrees of freedom, and symmetric, unless the symmetry is broken by
the elliptic equations.

This article is structured as follows. \Cref{sec:fluxform} details the generic
first-order flux formulation that serves as the starting point for our DG
discretization. \Cref{sec:dgdisc} develops the unified DG scheme. In
\cref{sec:tests} we apply the DG scheme to a set of increasingly challenging
test problems. The test problems include scenarios derived from general
relativity that feature sets of coupled, strongly nonlinear equations on a
curved manifold with nonlinear boundary conditions, solved on curved meshes. We
conclude in \cref{sec:conclusion}.

\section{First-order flux formulation}\label{sec:fluxform}

We consider second-order elliptic PDEs of one or more \enquote{primal} variables
$u_A(\bm{x})$, where the index $A$ labels the variables. The variables can be
scalars (like in the Poisson equation) or tensorial quantities (like in an
elasticity problem). We reduce the PDEs to first order by introducing
\enquote{auxiliary} variables $v_A(\bm{x})$, which typically are gradients of the
primal variables. We then restrict our attention to problems that can be
formulated in first-order flux form
\begin{equation}\label{eq:fluxform}
  -\partial_i \, \dgFi{\alpha}[u_A,v_A;\bm{x}] + \dgS_\alpha[u_A,v_A;\bm{x}] = \dgf_\alpha(\bm{x})
  \text{,}
\end{equation}
where the index $\alpha$ enumerates both $u_A$ and $v_A$. Here the fluxes
$\dgFi{\alpha}$ and the sources $\dgS_\alpha$ are functionals of the variables
$u_A$ and $v_A$, but not their derivatives, as well as the coordinates $\bm{x}$.
The fixed sources $\dgf_\alpha(\bm{x})$ are independent of the variables.
Lowercase Latin indices $i, j, k, l$ enumerate spatial dimensions, and we
employ the Einstein sum convention to sum over repeated indices.

The flux form~\eqref{eq:fluxform} is general enough to encompass a wide range of
elliptic problems. For example, a flat-space Poisson equation in Cartesian
coordinates
\begin{equation}\label{eq:poisson}
  -\partial_i \partial_i u(\bm{x}) = f(\bm{x})
\end{equation}
has the single primal variable $u(\bm{x})$. Choosing the auxiliary variable
$v_i = \partial_i u$ we can formulate the Poisson equation with the fluxes and
sources
\begin{subequations}\label{eq:poisson_fluxform}
\begin{alignat}{5}
  \label{eq:poisson_fluxform_aux}
  \tensor{\dgF}{_{\!v}^i_j} &= u \, \delta^i_j
  \text{,} &\quad
  \tensor{\dgS}{_v_{\,j}} &= v_j
  \text{,} &\quad
  \tensor{\dgf}{_v_{\,j}} &= 0
  \text{,} \\
  \label{eq:poisson_fluxform_prim}
  \dgFi{u} &= v_i
  \text{,} &\quad
  \dgS_u &= 0
  \text{,} &\quad
  \dgf_u &= f(\bm{x})
  \text{,}
\end{alignat}
\end{subequations}
where $\delta^i_j$ denotes the Kronecker delta. Note that
\cref{eq:poisson_fluxform_aux} is the definition of the auxiliary variable, and
\cref{eq:poisson_fluxform_prim} is the Poisson equation~\eqref{eq:poisson}.

The equation of linear elasticity in Cartesian coordinates,
\begin{equation}\label{eq:elasticity}
  -\partial_i \constrel^{ijkl} \partial_{(k}\,\displ_{l)} = f^j(\bm{x})
  \text{,}
\end{equation}
has the primal variable $\displ^i(\bm{x})$, describing the vectorial deformation
of an elastic material. The constitutive relation $\constrel^{ijkl}(\bm{x})$
captures the elastic properties of the material in the linear regime. Choosing
the symmetric
strain~$\strain_{ij}=\partial_{(i}\,\displ_{j)}=(\partial_i\,\displ_j +
\partial_j\,\displ_i) / 2$ as auxiliary variable we can formulate the elasticity
equation with the fluxes and sources
\begin{subequations}\label{eq:elasticity_fluxform}
\begin{alignat}{5}
  \label{eq:elasticity_fluxform_aux}
  \tensor{\dgF}{_{\!\strain}^i_j_k} &= \delta^i_{(j} \, \displ_{k)}
  \text{,} &\quad
  \tensor{\dgS}{_\strain_{\,j}_k} &= \strain_{jk}
  \text{,} &\quad
  \tensor{\dgf}{_\strain_{\,j}_k} &= 0
  \text{,} \\
  \label{eq:elasticity_fluxform_prim}
  \tensor{\dgF}{_{\!\displ}^i^j} &= \constrel^{ijkl} \strain_{kl}
  \text{,} &\quad
  \tensor{\dgS}{_\displ^j} &= 0
  \text{,} &\quad
  \tensor{\dgf}{_\displ^j} &= f^j(\bm{x})
  \text{.}
\end{alignat}
\end{subequations}
Again, \cref{eq:elasticity_fluxform_aux} is the definition of the auxiliary
variable and \cref{eq:elasticity_fluxform_prim} is the elasticity
equation~\eqref{eq:elasticity}. The fluxes and sources for the elasticity
system~\eqref{eq:elasticity_fluxform} have higher rank than those for the
Poisson system~\eqref{eq:poisson_fluxform}.

The first-order flux form~\eqref{eq:fluxform} also accommodates equations
formulated on a curved manifold which is equipped with a
metric~$g_{ij}(\bm{x})$. Such equations typically involve covariant
derivatives~$\nabla_i$ compatible with~$g_{ij}$. To formulate the equations in
flux form~\eqref{eq:fluxform} we expand covariant derivatives in partial
derivatives and Christoffel symbols $\Gamma^i_{jk} = \frac{1}{2} g^{il}
\left(\partial_j g_{kl} + \partial_k g_{jk} - \partial_l g_{jk}\right)$.
Christoffel symbols also appear when formulating equations in curvilinear
coordinates. In our scheme, the terms with partial derivatives are assigned to
the fluxes~$\dgF^i$ and the terms with Christoffel symbols are assigned to the
sources~$\dgS$. For example, a curved-space Poisson equation
\begin{equation}\label{eq:poisson_curved}
  -g^{ij}\nabla_i\nabla_j u(\bm{x}) = f(\bm{x})
\end{equation}
with auxiliary variable $v_i = \nabla_{\!i} u$ can be formulated with the
fluxes and sources
\begin{subequations}\label{eq:poisson_curved_fluxform}
\begin{alignat}{5}
  \tensor{\dgF}{_{\!v}^i_j} &= u \, \delta^i_j
  \text{,} &\quad
  \tensor{\dgS}{_v_{\,j}} &= v_j
  \text{,} &\quad
  \tensor{\dgf}{_v_{\,j}} &= 0
  \text{,} \\
  \dgFi{u} &= g^{ij} v_j
  \text{,} &\quad
  \dgS_u &= -\Gamma^i_{ij} g^{jk} v_k
  \text{,} &\quad
  \dgf_u &= f(\bm{x})
  \text{.}
\end{alignat}
\end{subequations}
Our strategy of expanding covariant derivatives differs from the formulations
employed for relativistic hyperbolic conservation laws
in~\ccite{Teukolsky2016-ja}, where fluxes are always vector fields and therefore
the covariant divergence can always be written in terms of partial derivatives
and the metric determinant.\footnote{See Eq.~(2.3) in \ccite{Teukolsky2016-ja}.}
In contrast, fluxes in the elliptic equations~\eqref{eq:fluxform} can be
higher-rank tensor fields, as exemplified in \cref{eq:elasticity_fluxform}.

The fixed sources~$\dgf_\alpha(\bm{x})$ could, in principle, be absorbed
in the sources~$\dgS_\alpha$. However, it is useful to keep these
variable-independent contributions separate for two reasons. First, they remain
constant throughout an elliptic solve, so they need not be recomputed when the
dynamic variables change. Second, the constant contributions represent a
nonlinearity in the variables~$u_A$ and~$v_A$ when included in the
sources~$\dgS_\alpha$. Assigning the constant contributions to the
fixed sources~$\dgf_\alpha$ eliminates this particular nonlinearity, hence
allowing us to avoid an explicit linearization procedure if the remaining
sources $\dgS_\alpha$ are linear.

The \hyperref[app:systems]{Appendix} lists fluxes and sources for selected
elliptic problems. Our focus on systems in generic first-order flux form
allows us to solve a variety of elliptic systems by only implementing their fluxes
and sources. We now proceed to discretize this generic formulation.

\section{DG discretization of the flux formulation}\label{sec:dgdisc}

In this section we develop the unified DG scheme for elliptic equations in flux
form, \cref{eq:fluxform}. Novel features of our scheme are the formulation of DG
residuals and boundary conditions in terms of generic fluxes and sources of
arbitrary tensor rank (\cref{sec:dgresiduals,sec:bc}), and the generalized
internal-penalty numerical flux (\cref{sec:numflux}). The Schur-complement
strategy of eliminating auxiliary degrees of freedom has been employed before,
e.g., in \ccite{Fortunato2019}, but we generalize it to a larger class of
equations, including equations with nonlinear fluxes or sources
(\cref{sec:compact}). We follow~\ccite{Teukolsky2016-ja}
whenever possible and refer to \ccite{HesthavenWarburton} for details that have
become standard in the DG literature.\footnote{\Ccite{Teukolsky2016-ja}
underpins the \emph{hyperbolic} DG formulations in the \spectre{} code.
Formulating elliptic and hyperbolic DG schemes in a similar way allows us to
share some of the DG implementation details.}

\subsection{Domain decomposition}\label{sec:domain}

We adopt the same domain decomposition based on deformed cubes detailed
in~\ccites{Teukolsky2016-ja,Kidder2017-nz,Vincent2019qpd} and summarize it here.

\begin{figure}
  \centering
  \subfloat[Domain\protect\label{fig:domain}]{
    \tikzinput{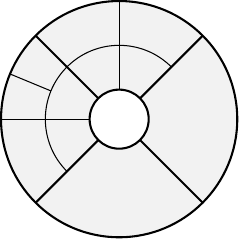}
  }\\
  \subfloat[Element\protect\label{fig:element}]{
    \tikzinput{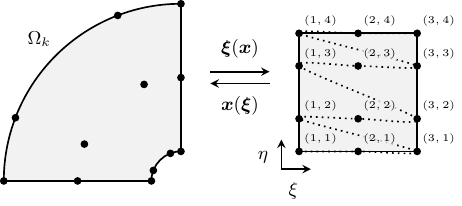}
  }
  \caption{
    \label{fig:domain_and_elements}
    \emph{Top:} Geometry of a two-dimensional computational domain composed of
    four wedge-shaped blocks. Each block is split in one or more nonoverlapping
    elements~$\Omega_k$. \emph{Bottom:} The coordinate
    transformation~$\bm{\xi}(\bm{x})$ maps an element to a reference
    cube~$[-1,1]^2$ with logical coordinate axes~$\bm{\xi}=(\xi,\eta)$. In this
    example we chose $N_{k,\xi}=3$ and $N_{k,\eta}=4$ Legendre-Gauss-Lobatto
    collocation points along $\xi$ and $\eta$, respectively. Each grid point is
    labeled with its index~$(p_\xi,p_\eta)$. The dotted line connects points in
    the order they are enumerated in by the index~$p$.}
\end{figure}

A $d$-dimensional computational domain $\Omega \subset \mathbb{R}^d$ is composed
of \emph{elements}~$\Omega_k \subset \Omega$ such that
$\Omega=\bigcup_k\Omega_k$. Elements do not overlap, but they share boundaries,
as illustrated in \cref{fig:domain}. Each element
carries an invertible map $\bm{\xi}(\bm{x})$ from the coordinates $\bm{x} \in
\Omega_k$, in which the elliptic equations~\eqref{eq:fluxform} are formulated,
to \emph{logical} coordinates $\bm{\xi} \in [-1, 1]^d$ representing a
$d$-dimensional reference cube. Inversely, $\bm{x}(\bm{\xi})$ maps the reference
cube to the element~$\Omega_k$. We define its Jacobian as
\begin{equation}
  \jac^i_j \defeq \pdv{{x^i}}{{\xi^j}}
\end{equation}
with determinant~$\jac$ and inverse~$\invjac^j_i = \partial \xi^j / \partial
x^i$.

Within each element~$\Omega_k$ we choose a set of $N_{k,i}$ grid points in every
dimension $i$. We place them at logical coordinates~$\xi_{p_i}$, where the index
$p_i\in\{1,\ldots,N_{k,i}\}$ identifies the grid point along dimension~$i$. The
points are laid out in a regular grid along the logical coordinate axes, so an
element has a total of $N_k=\prod_{i=1}^d N_{k,i}$ $d$-dimensional grid
points~$\bm{\xi}_p=(\xi_{p_1},\ldots,\xi_{p_d})$. The
index~$p\in\{1,\ldots,N_k\}$ identifies the grid point regardless of dimension.
The full domain has $N_\mathrm{points}=\sum_k N_k$ grid points. The grid points
within each element are not uniformly spaced in logical coordinates. Instead, we
choose Legendre-Gauss-Lobatto (LGL) collocation points, i.e., the
points~$\xi_{p_i}$ fall at the roots of the $(N_{k,i}-1)$th Legendre polynomial
plus a point on each side of the element, at~$-1$ and~$1$.\footnote{See, e.g.,
Algorithm~25 in \ccite{Kopriva}.} It is equally possible to choose
Legendre-Gauss (LG) collocation points, i.e., the roots of the $N_{k,i}$th
Legendre polynomial.\footnote{See, e.g., Algorithm~23 in \ccite{Kopriva}.}
\Cref{fig:element} illustrates the geometry of an element.

Fields are represented numerically by their values at the grid points. To
facilitate this we construct the one-dimensional Lagrange polynomials
\begin{equation}\label{eq:lagr}
  \lagr_{p_i}(\xi)\defeq \prod_{\substack{q_i=1 \\ q_i\neq p_i}}^{N_{k,i}}
  \frac{\xi - \xi_{q_i}}{\xi_{p_i} - \xi_{q_i}}
  \quad \text{with} \quad \xi \in [-1,1]
\end{equation}
and employ their product to define the $d$-dimensional basis functions
\begin{equation}\label{eq:lagrprod}
  \bas_p(\bm{\xi}) \defeq \prod_{i=1}^d \lagr_{p_i}(\xi^i)
  \quad \text{with} \quad \bm{\xi} \in [-1,1]^d
  \text{.}
\end{equation}
The choice of Lagrange polynomials makes \cref{eq:lagrprod} a nodal basis with
the useful property~$\bas_p(\bm{\xi}_q)=\delta_{pq}$. We use the nodal
basis~\eqref{eq:lagrprod} to approximate any field $u(\bm{x})$ within an
element~$\Omega_k$ by its discretization
\begin{equation}\label{eq:field_expansion}
  u^{(k)}(\bm{x}) \defeq \sum_{p=1}^{N_k} u_p \bas_p(\bm{\xi}(\bm{x}))
  \quad \text{with} \quad \bm{x} \in \Omega_k
  \text{,}
\end{equation}
where the coefficients $u_p=u(\bm{x}(\bm{\xi}_p))$ are the field values at
the grid points. We denote the set of discrete field values within an element
$\Omega_k$ as
\begin{equation}
  \grd{u}^{(k)}=(u_1,\ldots,u_{N_k})
  \text{,}
\end{equation}
and the collection of discrete field values over \emph{all} elements
as~$\grd{u}$. The
discretization~\eqref{eq:field_expansion} approximates fields with polynomials
of degree~$(N_{k,i}-1)$ in dimension~$i$. Although rarely needed, field values
at other points within an element can be obtained by Lagrange
interpolation~\eqref{eq:field_expansion}. The field values at element
boundaries are double valued because the Lagrange interpolation from
neighboring elements to their shared boundary is double valued. Therefore,
field approximations will in general be discontinuous at element boundaries.

The test problems in \cref{sec:tests} illustrate a few examples of domain
decompositions. We refer the reader to, e.g., \ccite{HesthavenWarburton} for
further details on the choice of collocation points, basis functions and their
relation to spectral properties of DG schemes.

\subsection{DG residuals}\label{sec:dgresiduals}

The DG residuals represent the set of equations to be solved for the discrete
primal field values $\grd{u}_A$. The derivation in this section follows
the standard procedure, e.g., laid out in \ccite{HesthavenWarburton}, applied to
the generic elliptic flux formulation~\eqref{eq:fluxform}, and taking details
such as a curved manifold into account.

In the spirit of a Galerkin scheme we project our target
PDEs~\eqref{eq:fluxform} onto the same set of basis functions $\bas_p(\bm{\xi})$
that is used to approximate fields within an element~$\Omega_k$,
\begin{equation}\label{eq:dgres_start}
  -(\bas_p, \partial_i \dgF^i)_{\Omega_k} + (\bas_p, \dgS)_{\Omega_k} = (\bas_p, \dgf)_{\Omega_k}
  \text{.}
\end{equation}
Here we dropped the index $\alpha$ that enumerates the equations, and we define
the inner product on~$\Omega_k$,
\begin{subequations}\label{eq:dgproj}
\begin{align}
  (\phi, \pi)_{\Omega_k} \defeq{}& \int_{\Omega_k} \phi(\bm{x}) \pi(\bm{x}) \sqrt{g} \dd{^dx} \\
  ={}& \int_{[-1,1]^d} \phi(\bm{x}(\bm{\xi})) \pi(\bm{x}(\bm{\xi})) \sqrt{g} \, \jac \dd{^d\xi}
  \text{.}
\end{align}
\end{subequations}
These integrals are defined with respect to proper volume
$\dd{V}=\sqrt{g}\dd{^dx}=\sqrt{g}\,\jac\dd{^d\xi}$, where $g$ denotes the
metric determinant in the coordinates $\bm{x}$ in which
\cref{eq:fluxform} is formulated. It refers to the metric that
covariant derivatives in the equations are compatible with. Since the basis
polynomials, \cref{eq:lagrprod}, are functions of logical coordinates, we
abbreviate $\bas_p(\bm{\xi}(\bm{x}))$ with $\bas_p(\bm{x})$ here.

The terms without derivatives in \cref{eq:dgres_start} are straightforward to
discretize. We approximate the field~$\dgf$, or similarly~$\dgS$, using the
expansion in basis functions~\eqref{eq:field_expansion} to find
\begin{equation}
  (\bas_p, \dgf)_{\Omega_k} \approx (\bas_p, \bas_q)_{\Omega_k} \dgf_q = \dgM_{pq} \dgf_q
  \text{,}
\end{equation}
using the symmetric \emph{mass matrix} on the element~$\Omega_k$,
\begin{subequations}\label{eq:massmat}
\begin{align}
  \dgM_{pq} \defeq{}& (\bas_p, \bas_q)_{\Omega_k} \\
  ={}& \int_{[-1,1]^d} \bas_p(\bm{\xi}) \bas_q(\bm{\xi}) \sqrt{g} \, \jac \dd{^d\xi}
  \text{.}
\end{align}
\end{subequations}
We will discuss strategies to evaluate the mass matrix on the elements of the
computational domain in \cref{sec:linop_impl}.

The divergence term in \cref{eq:dgres_start} encodes the principal part of the
elliptic PDEs and requires more care in its discretization. The derivatives in
this term will help us couple grid points across element boundaries. To this end
we integrate by parts to obtain a boundary term
\begin{equation}\label{eq:weakdiv_raw}
  (\bas_p, \partial_i \dgF^i)_{\Omega_k} = -(\partial_i \bas_p, \dgF^i)_{\Omega_k}
  + (\bas_p, n_i \dgF^i)_{\partial\Omega_k}
  \text{,}
\end{equation}
where $n_i$ is the outward-pointing unit normal one form on the element boundary
$\partial\Omega_k$. The \emph{unnormalized} face normal is computed from the
Jacobian as $\tilde{n}_i=\mathrm{sgn}(\xi^j)\invjac^j_i$, where $\xi^j$ is the
logical coordinate that is constant on the particular face and no sum over~$j$
is implied. The face normal is normalized as $n_i = \tilde{n}_i /
\sqrt{\tilde{n}_k \tilde{n}_l g^{kl}}$ using the inverse metric
$g^{ij}(\bm{x})$. The surface integral in \cref{eq:weakdiv_raw} is defined just
like \cref{eq:dgproj},
\begin{subequations}
\begin{align}
  (\phi, \pi)_{\partial\Omega_k} \defeq{}& \int_{\partial\Omega_k} \phi(\bm{x}) \pi(\bm{x}) \sqrt{\surf{g}} \dd{^{d-1}x} \\
  ={}& \int_{[-1,1]^{d-1}} \phi(\bm{x}(\bm{\xi})) \pi(\bm{x}(\bm{\xi})) \sqrt{\surf{g}} \, \surf{\jac} \dd{^{d-1}\xi}
  \text{,}
\end{align}
\end{subequations}
using the element boundary's $(d-1)$-dimensional proper volume
$\dd{\Sigma}=\sqrt{\surf{g}}\dd{^{d-1}x}=\sqrt{\surf{g}}\,\surf{\jac}\dd{^{d-1}\xi}$,
where $\surf{g}$ is the surface metric determinant induced by the
metric~$g_{ij}$ and $\surf{\jac}$ is the surface Jacobian.

The crucial step that couples grid points across element boundaries follows from
the field $n_i\dgF^i$ being double valued on any section of the boundary that an
element shares with a neighbor, with one value arising from either side. We must
make a choice how to combine the two values
from either side of a shared element boundary. This choice is often referred to
as a \emph{numerical flux}. For now we will denote the function that combines
values from both sides of a boundary as $(n_i\dgF^i)^*$ and refer to
\cref{sec:numflux} for details on our particular choice of numerical flux.
Substituting the numerical flux in \cref{eq:weakdiv_raw} yields the \emph{weak}
form of the equations,
\begin{equation}\label{eq:weakdiv}
  (\bas_p, \partial_i \dgF^i)_{\Omega_k} = -(\partial_i \bas_p, \dgF^i)_{\Omega_k}
  + (\bas_p, (n_i \dgF^i)^*)_{\partial\Omega_k}
  \text{.}
\end{equation}
The numerical flux in \cref{eq:weakdiv} introduces a coupling between
neighboring elements that allows us to obtain numerical solutions spanning the
full computational domain. Another integration by parts of \cref{eq:weakdiv}
yields the \emph{strong} form of the equations,
\begin{equation}\label{eq:strongdiv}
  (\bas_p, \partial_i \dgF^i)_{\Omega_k} = (\bas_p, \partial_i\dgF^i)_{\Omega_k}
  + (\bas_p, (n_i \dgF^i)^* - n_i \dgF^i)_{\partial\Omega_k}
  \text{.}
\end{equation}
We will make use of both the strong and the weak form to obtain symmetric DG
operators (see \cref{sec:symm}). Approximating~$\dgF^i$ using its expansion
in basis functions~\eqref{eq:field_expansion} we find
\begin{equation}
  (\bas_p, \partial_i \dgF^i)_{\Omega_k} \approx
  (\bas_p, \partial_i\bas_q)_{\Omega_k} \dgF^i_q =
  \dgMD_{i,pq} \dgF^i_q
  \text{,}
\end{equation}
where the \emph{stiffness matrix} on the element~$\Omega_k$ is
\begin{subequations}\label{eq:stiffmat}
\begin{align}
  \dgMD_{i,pq} \defeq{}& (\bas_p, \partial_i\bas_q)_{\Omega_k} \\
  ={}& \int_{[-1,1]^d} \bas_p(\bm{\xi}) \, \pdv{\bas_q}{{\xi^j}}\!(\bm{\xi})
  \, \invjac^j_i \sqrt{g} \, \jac \dd{^d\xi}
  \text{.}
\end{align}
\end{subequations}
The divergence term in its weak form can be expressed in terms of the
stiffness-matrix transpose $\dgMD^T_{i,pq}=\dgMD_{i,qp}$,
\begin{equation}
  -(\partial_i\bas_p, \dgF^i)_{\Omega_k} \approx
  -(\partial_i\bas_p, \bas_q)_{\Omega_k} \dgF^i_q =
  -\dgMD^T_{i,pq} \dgF^i_q
  \text{.}
\end{equation}
Evaluation of the stiffness matrix and its transpose is discussed in
\cref{sec:linop_impl}.

We now turn towards discretizing the last remaining piece of the DG residuals,
the boundary integrals in \cref{eq:weakdiv,eq:strongdiv}. It
involves a \enquote{lifting} operation: the integral only depends on field values on the
element boundary but it may contribute to every component~$p$ of the DG
residual, hence it is \enquote{lifted} to the volume. However, on an LGL grid all
components~$p$ that correspond to grid points away from the boundary evaluate to
zero because they contain at least one Lagrange polynomial that vanishes at the
boundary collocation point. This is not the case on an LG grid, where evaluating
the Lagrange polynomials on the boundary produces an interpolation into the
volume. Expanding the boundary fluxes in basis
functions~\eqref{eq:field_expansion} we find
\begin{equation}
  (\bas_p, n_i\dgF^i)_{\partial\Omega_k} \approx
  (\bas_p, \bas_q)_{\partial\Omega_k} (n_i\dgF^i)_q =\dgML_{pq}(n_i\dgF^i)_q
  \text{,}
\end{equation}
where we have defined the \emph{lifting operator} on the element~$\Omega_k$,
\begin{subequations}\label{eq:liftop}
\begin{align}
  \dgML_{pq} \defeq{}& (\bas_p,\bas_q)_{\partial\Omega_k} \\
  ={}& \int_{[-1,1]^{d-1}} \bas_p(\bm{\xi})\bas_q(\bm{\xi})
  \sqrt{\surf{g}} \, \surf{\jac} \dd{^{d-1}\xi}
  \text{.}
\end{align}
\end{subequations}
\Cref{sec:linop_impl} provides details on evaluating the lifting operator.

Assembling the pieces of the discretization and restoring the index $\alpha$
that enumerates the equations, the DG residuals on the element~$\Omega_k$ 
in strong form are
\begin{subequations}\label{eq:dg_op_full}
\begin{align}
  -\dgMD_i \cdot \dgFi{\alpha}
  - \dgML \cdot ((n_i \dgFi{\alpha})^*& - n_i \dgFi{\alpha}) \nonumber \\
  &+ \dgM \cdot \dgS_\alpha
  = \dgM \cdot \dgf_\alpha \label{eq:dg_op_full_strong}
  \text{,}
\end{align}
where~$\cdot$ denotes a matrix multiplication with the field values over the
computational grid of an element. The DG residuals in weak form are
\begin{equation}
  \dgMD^T_i \cdot \dgFi{\alpha}
  - \dgML \cdot (n_i \dgFi{\alpha})^*
  + \dgM \cdot \dgS_\alpha
  = \dgM \cdot \dgf_\alpha \label{eq:dg_op_full_weak}
  \text{.}
\end{equation}
\end{subequations}
We can choose either the strong or the weak form for each variable~$\alpha$.

Since the fluxes and
sources are computed from the primal and auxiliary variables, the DG
residuals~\eqref{eq:dg_op_full} are algebraic equations for the discrete values
$\grd{u}_A$ and~$\grd{v}_A$ on all elements and grid points in the
computational domain. The left-hand side of \cref{eq:dg_op_full} is an
operator $\dgA(\grd{u}_A, \grd{v}_A)$ and the right-hand side of
\cref{eq:dg_op_full} is a fixed value at every grid point, so
\cref{eq:dg_op_full} has the structure
\begin{equation}\label{eq:dg_res_full}
  \dgA(\grd{u}_A, \grd{v}_A) = \grd{\dgb}
  \text{.}
\end{equation}
If the fluxes and sources are linear, the DG
operator $\dgA(\grd{u}_A, \grd{v}_A)$ can be represented as a square
matrix, and \cref{eq:dg_res_full} is a matrix equation. The size of
the DG operator $\dgA(\grd{u}_A, \grd{v}_A)$ is the product
of~$N_\mathrm{points}$ with the number of both primal and auxiliary variables.

\begin{figure}
  \centering
  \includegraphics[width=\columnwidth]{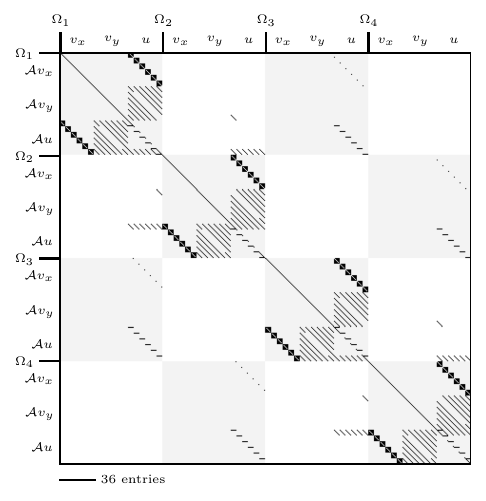}
  \caption{
  \label{fig:noncompact_operator}
  Matrix representation of the noncompact DG operator in strong
  form~\eqref{eq:dg_op_full_strong} for a two-dimensional Poisson
  equation~\eqref{eq:poisson_fluxform}. The computational domain is partitioned
  into $2 \times 2$ elements with $6 \times 6$ LGL grid points each. The image
  shows the nonzero entries of the operator matrix, i.e., its sparsity pattern,
  in the order laid out in~\cref{fig:element}.}
\end{figure}

\Cref{fig:noncompact_operator} presents a visualization of the DG
operator $\dgA(\grd{u}_A, \grd{v}_A)$ for a Poisson equation on a
regular grid. The axes annotate entries of the operator that correspond to the
\enquote{input} variables $v_i$ and $u$, and to the corresponding \enquote{output} DG residuals.
The mass matrix applied to $v_i$ appears as a diagonal line (see
\cref{sec:linop_impl}) and the stiffness matrices applied to both $v_i$ and $u$
appear as block-diagonal and shaded regions for derivatives in $x$ and $y$,
respectively. The remaining entries represent the coupling between neighboring
elements through the numerical flux (see \cref{sec:numflux}). Note that the
elements $\Omega_1$ and $\Omega_4$ as well as $\Omega_2$ and $\Omega_3$
decouple, because they share no boundaries as they are placed diagonally across
the $2 \times 2$ grid of elements. Solving the Poisson equation amounts to
inverting the matrix pictured in \cref{fig:noncompact_operator}. However, it is
significantly cheaper to invert the equivalent compact operator pictured in
\cref{fig:compact_operator}, which we derive in the following section.

\subsection{Eliminating auxiliary degrees of freedom}\label{sec:compact}

So far we have treated the primal and the auxiliary equations of the
first-order formulation on the same footing, which means the discretized DG
operator applies to the primal variables as well as to the auxiliary variables.
However, the auxiliary equations inflate the size of the operator significantly,
increasing both its memory usage and the computational cost for solving it. In
this section we eliminate the auxiliary degrees of freedom from the DG operator,
demoting them to quantities that are only computed temporarily.

Many publications on DG formulations adopt a \enquote{primal
formulation} to eliminate auxiliary degrees of freedom from the DG
operator.\footnote{See, e.g., Section~7.2.2 in \ccite{HesthavenWarburton} or
Section~3 in \ccite{Arnold2002} for derivations of primal formulations for the
Poisson equation.} However, in practice we have found a simpler approach taking
a Schur complement of the discretized equations in flux form, e.g., applied in
\ccite{Fortunato2019}, more suited to the generic implementation of DG schemes.
The resulting DG operator remains equivalent to the original operator; i.e., it
has the same solutions up to numerical precision. This strategy is facilitated
by the auxiliary equations defining the auxiliary variables~$v_A$, such as
\cref{eq:poisson_fluxform_aux,eq:elasticity_fluxform_aux}. We assume
here that the auxiliary fluxes depend only on the primal variables,
$\dgFi{v_A}=\dgFi{v_A}[u_A;\bm{x}]$, and the auxiliary sources have the form
\begin{equation}\label{eq:auxsource}
  \dgS_{v_A}=v_A + \tilde{\dgS}_{v_A}[u_A;\bm{x}]
  \text{,}
\end{equation}
where $\tilde{\dgS}_{v_A}$ depends only on the primal variables. We further
assume $\dgf_{v_A}=0$ for convenience. All elliptic systems that we consider in
this article fulfill these assumptions. We insert \cref{eq:auxsource} into the
strong DG residuals~\eqref{eq:dg_op_full_strong} and solve for $v_A$ by
inverting the mass matrix to find
\begin{equation}\label{eq:aux_matform}
  v_A = \dgD_i \cdot \dgFi{v_A}
  + \dgL \cdot ((n_i\dgFi{v_A})^* - n_i\dgFi{v_A})
  - \tilde{\dgS}_{v_A}
  \text{,}
\end{equation}
where we define $\dgD_i \defeq \dgM^{-1}\dgMD_i$ and $\dgL \defeq
\dgM^{-1}\dgML$. Note that the right-hand side of \cref{eq:aux_matform} depends
only on the primal variables $u_A$. Evaluating the DG residuals now amounts to
first computing the auxiliary variables $v_A$ by \cref{eq:aux_matform}, and then
using them to evaluate the primal equations. This approach preserves the freedom
to use either the strong form~\eqref{eq:dg_op_full_strong} for the primal
equations,
\begin{subequations}\label{eq:dg_op_compact}
\begin{align}
  -\dgMD_i \cdot \dgFi{u_A}
  - \dgML \cdot ((n_i \dgFi{u_A})^* &- n_i \dgFi{u_A}) \nonumber \\
  &+ \dgM  \cdot \dgS_{u_A}
  = \dgM \cdot \dgf_{u_A} \label{eq:dg_op_compact_strong}
  \text{,}
\end{align}
or the weak form~\eqref{eq:dg_op_full_weak},
\begin{equation}
  \dgMD^T_i \cdot \dgFi{u_A}
  - \dgML \cdot (n_i \dgFi{u_A})^* 
  + \dgM \cdot \dgS_{u_A} 
  = \dgM \cdot \dgf_{u_A} \label{eq:dg_op_compact_weak}
  \text{.}
\end{equation}
\end{subequations}
The strong scheme is slightly easier to implement because the primal and
auxiliary equations involve the same set of operators. The strong-weak scheme,
i.e., selecting the strong form for the auxiliary
equations~\eqref{eq:aux_matform} and the weak form for the primal
equations~\eqref{eq:dg_op_compact_weak}, has the advantage that the DG operator
can be symmetric as discussed in \cref{sec:symm}.

\begin{figure}
  \centering
  \includegraphics[width=\columnwidth]{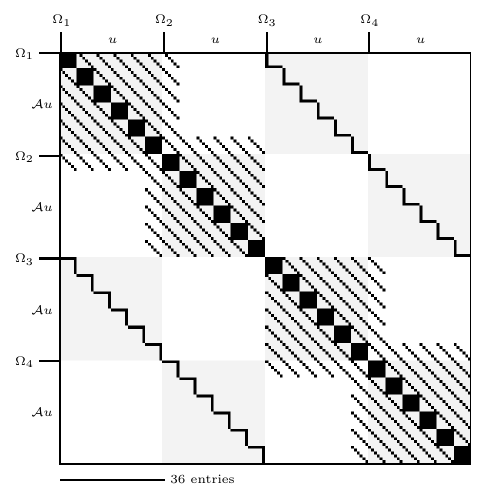}
  \caption{
  \label{fig:compact_operator}
  Matrix representation of the compact DG operator in strong
  form~\eqref{eq:dg_op_compact_strong} for the two-dimensional Poisson problem
  detailed in \cref{fig:noncompact_operator}. No auxiliary degrees of freedom
  inflate the size of the operator. This matrix is the Schur complement to the
  matrix pictured in \cref{fig:noncompact_operator}.}
\end{figure}

\begin{figure}
  \centering
  \includegraphics[width=\columnwidth]{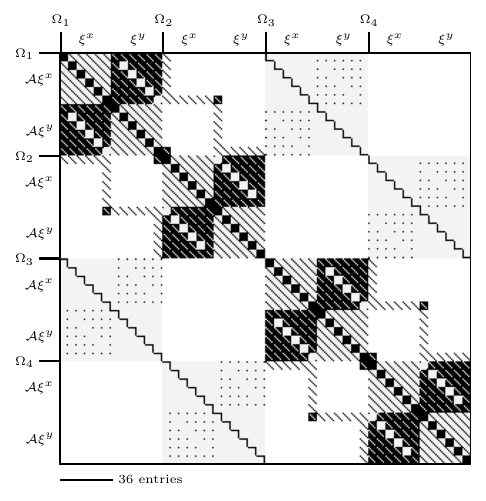}
  \caption{
  \label{fig:elasticity_operator}
  Matrix representation of the compact DG operator in strong
  form~\eqref{eq:dg_op_compact_strong} for a two-dimensional elasticity
  problem~\eqref{eq:elasticity_fluxform} with an isotropic-homogeneous
  constitutive relation~$\constrel^{ijkl}$. The computational domain is again
  partitioned into $2 \times 2$ elements with $6 \times 6$ LGL grid points
  each.}
\end{figure}

For linear equations the strategy employed in \cref{eq:aux_matform} of
eliminating the auxiliary variables is equivalent to taking a Schur complement
of the DG operator with respect to the (invertible) mass matrix, but the
strategy works for nonlinear equations as well. The result is an
operator $\dgA(\grd{u}_A)$ of only the discrete primal variables on all
elements and grid points in the computational domain, so the DG
residuals~\eqref{eq:dg_op_compact} have the form
\begin{equation}
  \dgA(\grd{u}_A) = \grd{b}
  \text{.}
\end{equation}
The size of the DG operator $\dgA(\grd{u}_A)$
is the product of $N_\mathrm{points}$ with the number of primal variables. No
auxiliary degrees of freedom from the first-order formulation inflate the size
of the operator. We refer to such DG operators~$\dgA(\grd{u}_A)$ of only
the primal degrees of freedom as \emph{compact} if they also only involve
couplings between nearest-neighbor elements~\cite{Peraire2007-ff}. The coupling
between elements is related to the choice of numerical flux~$(n_i
\dgFi{\alpha})^*$ and the subject of \cref{sec:numflux}. If the fluxes and
sources are linear, $\dgA(\grd{u}_A)$ can be represented as a square
matrix.

\Cref{fig:compact_operator,fig:elasticity_operator} present visualizations
of~$\dgA$ for a Poisson equation and an elasticity equation, respectively.
The block-diagonal structure in~\cref{fig:compact_operator} represents the
DG-discretized two-dimensional Laplacian on the four elements of the
computational domain. The entries that break the block-diagonal structure
represent the coupling between nearest-neighbor elements through the numerical
flux (\cref{sec:numflux}).

\subsection{A generalized internal-penalty numerical flux}\label{sec:numflux}

Up to this point we have not made a choice for the numerical flux $(n_i
\dgF^i)^*$ that combines double-valued field values on element boundaries. The
numerical flux is a function of the field values on both sides of the boundary.
From the perspective of one of the two adjacent elements we refer to the field
values on itself as \emph{interior} and to the field values the neighboring
element as \emph{exterior}. Contrary to much of the DG literature we formulate
the numerical flux entirely in terms of the primal and auxiliary \emph{boundary
flux} quantities $n_i \dgFi{u_A}$ and $n_i \dgFi{v_A}$ on either side of the
boundary instead of the primal and auxiliary variables $u_A$ and $v_A$. This
choice keeps our scheme applicable to the wide range of elliptic problems
defined in \cref{sec:fluxform}. The numerical flux presented here is a
generalization of the symmetric internal penalty (SIP) scheme that is widely
used in the literature \cite{Douglas1976-hx, Arnold2002, HesthavenWarburton,
Vincent2019qpd}. Our \emph{generalized internal-penalty numerical flux} is
\begin{subequations}\label{eq:numflux}
\begin{align}
  (n_i \dgFi{v_A})^* = \frac{1}{2}\Bigl[
    &n_i^\mathrm{int} \dgFi{v_A}(u_A^\mathrm{int}) -
    n_i^\mathrm{ext} \dgFi{v_A}(u_A^\mathrm{ext})\Bigr]
    \label{eq:numfluxaux}\text{,}\\
  (n_i \dgFi{u_A})^* =
    \frac{1}{2}\Bigl[
    &n_i^\mathrm{int} \dgFi{u_A}\bigl(\partial_j \dgFj{v_A}(u_A^\mathrm{int}) - \tilde{\dgS}_{v_A}(u_A^\mathrm{int})\bigr) \nonumber \\
     &- n_i^\mathrm{ext} \dgFi{u_A}\bigl(\partial_j \dgFj{v_A}(u_A^\mathrm{ext}) - \tilde{\dgS}_{v_A}(u_A^\mathrm{ext})\bigr)\Bigr] \nonumber \\
    - \sigma \Bigl[
      &n_i^\mathrm{int} \dgFi{u_A}\bigl(n_j^\mathrm{int} \dgFj{v_A}(u_A^\mathrm{int})\bigr) \nonumber \\
       &-n_i^\mathrm{ext} \dgFi{u_A}\bigl(n_j^\mathrm{ext} \dgFj{v_A}(u_A^\mathrm{ext})\bigr)\Bigr] \label{eq:numfluxprimal}
  \text{.}
\end{align}
\end{subequations}
The numerical flux for the auxiliary equations, \cref{eq:numfluxaux}, averages
the boundary fluxes of the two adjacent elements. The numerical flux for the
primal equations, \cref{eq:numfluxprimal}, is an average augmented with a
penalty contribution with parameter $\sigma$.

Note that the numerical flux~\eqref{eq:numflux} involves only the primal fields
and their derivatives, and thus is independent of the auxiliary fields
altogether, as is typical for internal-penalty schemes. This has the practical
advantage that the contribution from either side of the boundary to both the
primal and the auxiliary numerical flux in
\cref{eq:aux_matform,eq:dg_op_compact} can be computed early in the algorithm
and communicated together, coupling only nearest-neighbor elements and thus
making the DG operator compact. If the primal numerical
flux~\eqref{eq:numfluxprimal} depended on the auxiliary fields, evaluating the
DG operator~\eqref{eq:dg_op_compact} would require a separate communication once
the boundary corrections have been added to the auxiliary
equation~\eqref{eq:aux_matform}, effectively coupling nearest-neighbor elements
as well as next-to-nearest-neighbor elements.\footnote{Couplings to
next-to-nearest-neighbor elements is a well-known disadvantage of LDG-type
(\enquote{local discontinuous Galerkin}) numerical fluxes and has led to the
development of compact schemes such as~\ccite{Peraire2007-ff}.}

DG literature usually assumes that the face normals on either side of the
boundary are exactly opposite, $n_i^\mathrm{ext} = -n_i^\mathrm{int}$. This
assumption breaks when the background geometry responsible for the normalization
of face normals depends on the dynamic variables, since those are
discontinuous across the boundary. All of the elliptic problems that we are
expecting to solve in the near future are formulated on a fixed background
geometry, but it is useful to distinguish between the interior and exterior face
normals nonetheless because the quantity $n_i\dgF^i$ is cheaper to communicate
than $\dgF^i$. Therefore, we always project an element's boundary fluxes onto
the face normal before communicating the quantity.

For a simple flat-space Poisson system~\eqref{eq:poisson_fluxform} our
generalized internal-penalty numerical flux~\eqref{eq:numflux} reduces to the
canonical SIP,\footnote{See Eq.~(3.21) in~\ccite{Arnold2002} or Section~7.2
in~\ccite{HesthavenWarburton}.}
\begin{subequations}
\begin{align}
  (n_i \dgFi{v_j})^* &= n_j^\mathrm{int} u^* = \frac{1}{2} n_j^\mathrm{int} \big(
    u^\mathrm{int} + u^\mathrm{ext}\big) \\
  (n_i \dgFi{u})^* &= n_i^\mathrm{int} {v^i}^* = \begin{aligned}[t]
    &\frac{1}{2} n_i^\mathrm{int} \big(\partial^i u^\mathrm{int} + \partial^i u^\mathrm{ext}\big) \\
    &- \sigma \big(u^\mathrm{int} - u^\mathrm{ext}\big)
    \text{.}
  \end{aligned}
\end{align}
\end{subequations}

As is well studied for the canonical SIP numerical flux, the penalty
factor~$\sigma$ is responsible for removing zero eigenmodes and impacts the
conditioning of the linear operator to be
solved~\cite{HesthavenWarburton,Shahbazi2005-dr}. It scales inversely with the
element size~$h$ and quadratically with the polynomial degree~$p$, both
orthogonal to the element face.\footnote{See~\ccite{Shahbazi2005-dr} for sharp
results for the optimal penalty factor on triangular and tetrahedral meshes, and
Table~3.1 in~\ccite{Hillewaert2013-jn} for a generalization to hexahedral
meshes.} Both~$h$ and~$p$ can be different on either side of the element
boundary, so we choose
\begin{equation}\label{eq:penalty}
  \sigma = \penaltyparam \, \frac{\big(\max(p^\mathrm{int},p^\mathrm{ext})+1\big)^2}{\min(h^\mathrm{int}, h^\mathrm{ext})}
  \text{,}
\end{equation}
where we follow~\ccite{Hillewaert2013-jn} in choosing the scaling with the
polynomial degree~$p$ on hexahedral meshes, and we follow~\ccite{Vincent2019qpd}
in our definition of the element size~$h=2 / \sqrt{\tilde{n}_i\tilde{n}_j
g^{ij}}$.\footnote{Note that~\ccite{Vincent2019qpd} omits the factor of $2$ in
the definition of~$h$, which we include so the definition reduces to the
canonical element size on rectilinear meshes.} Note that~$h$ generally varies
over the element face on curved meshes or on a curved manifold, and that the
$\min$ operation in \cref{eq:penalty} is taken pointwise, so
$\sigma$ also varies over the element face. The remaining \emph{penalty
parameter}~$\penaltyparam \geq 1$ remains freely specifiable. Note also that we
do not need to include a problem-specific scale in the penalty factor, as is
done in \ccites{Hartmann2008-zi,Epshteyn2007-dm,De_Basabe2008-wf}, because the
generic numerical flux~\eqref{eq:numfluxprimal} already includes such scales in
the fluxes~$\dgF^i$.

\subsection{Imposing boundary conditions}\label{sec:bc}

The flux formulation allows imposing a wide range of boundary conditions
relatively easily \enquote{through the fluxes} without the need to treat external
boundaries any differently than internal boundaries between neighboring
elements. Imposing boundary conditions amounts to specifying the exterior
quantities in the numerical flux, \cref{eq:numflux}. This strategy is often
referred to as imposing boundary conditions through \enquote{ghost} elements. As
suggested in, e.g., \ccite{HesthavenWarburton}, we impose boundary conditions on
the \emph{average} of the boundary fluxes to obtain faster convergence.
Therefore, on external boundaries, we choose for the exterior quantities in the
numerical flux~\eqref{eq:numflux}
\begin{equation}\label{eq:bc}
  (n_i\dgFi{\alpha})^\mathrm{ext} = (n_i\dgFi{\alpha})^\mathrm{int} - 2 (n_i\dgFi{\alpha})^\mathrm{b}
  \text{,}
\end{equation}
where we set the quantities~$(n_i\dgFi{\alpha})^\mathrm{b}$ according to the
boundary conditions at hand. Here we define $n_i^\mathrm{b}=n_i^\mathrm{int}$,
i.e., we choose to compute external boundary fluxes with the face normal
pointing \emph{out} of the computational domain. The symmetry between the primal
and the auxiliary equations in the first-order flux
formulation~\eqref{eq:fluxform} that we employ throughout this article makes
this approach of imposing boundary conditions particularly straightforward: a
choice of \emph{auxiliary} boundary fluxes $(n_i\dgFi{v_A})^\mathrm{b}$ imposes
Dirichlet-type boundary conditions and a choice of \emph{primal} boundary fluxes
$(n_i\dgFi{u_A})^\mathrm{b}$ imposes Neumann-type boundary conditions. The
choice between Dirichlet-type and Neumann-type boundary conditions can be made
separately for every primal variable $u_A$ and every external boundary face,
simply by setting either $(n_i\dgFi{u_A})^\mathrm{b}$ or
$(n_i\dgFi{v_A})^\mathrm{b}$ and setting the remaining boundary fluxes to their
interior values $(n_i\dgF^i)^\mathrm{b} = (n_i\dgF^i)^\mathrm{int}$. Note that
we neither need to distinguish between primal and auxiliary variables in
\cref{eq:bc}, nor take the choice of Dirichlet-type or Neumann-type boundary
conditions into account, but we require only that $(n_i\dgF^i)^\mathrm{b}$ be chosen
appropriately for every variable. Then, the Neumann-type boundary conditions
enter the DG residuals directly through the numerical flux in
\cref{eq:dg_op_compact}, and the Dirichlet-type boundary conditions enter the DG
residuals through the numerical flux in \cref{eq:aux_matform}, which is
substituted in \cref{eq:dg_op_compact}. 

In practice, this setup means we can initialize $(n_i\dgFi{\alpha})^\mathrm{b} =
(n_i\dgFi{\alpha})^\mathrm{int}$ for all variables on a particular external
boundary face when preparing to apply the numerical flux, decide which boundary
fluxes to modify based on the boundary conditions we wish to impose on the
particular face, and then evaluate \cref{eq:bc} to compute the exterior
quantities in the numerical flux~\eqref{eq:numflux}. To impose Neumann-type
boundary conditions we set the primal boundary fluxes
$(n_i\dgFi{u_A})^\mathrm{b}$ directly, but to impose Dirichlet-type boundary
conditions we typically choose the primal field values $u_A^\mathrm{b}$ and
compute the auxiliary boundary fluxes as
$(n_i\dgFi{v_A})^\mathrm{b}=n_i^\mathrm{int} \dgFi{v_A}(u_A^\mathrm{b})$.

The auxiliary (Dirichlet-type) external boundary fluxes may depend on the
interior primal fields $u_A^\mathrm{int}$, and the primal (Neumann-type)
external boundary fluxes may depend on both the interior primal fields
$u_A^\mathrm{int}$ as well as the interior auxiliary fields $v_A^\mathrm{int}$.
This means we can impose a wide range of boundary conditions that may depend
linearly or nonlinearly on the dynamic fields. For example, a Robin boundary
condition for the Poisson equation~\eqref{eq:poisson}
or~\eqref{eq:poisson_curved},
\begin{equation}
  a \, u + b \, n^i \partial_i u = g(\bm{x}) \quad \text{on $\partial \Omega$}
  \text{,}
\end{equation}
where $a$ and $b$ are constants and $g(\bm{x})$ is a function defined on the
boundary, can be implemented as Neumann-type for $b \neq 0$,
\begin{equation}
  (n_i\dgFi{u})^\mathrm{b}=\frac{1}{b} \big(g(\bm{x}) - a \, u^\mathrm{int}\big)
  \text{,}
\end{equation}
and as Dirichlet-type for $b = 0$,
\begin{equation}
  (n_i\dgFi{v})^\mathrm{b} = n_i^\mathrm{int} \dgFi{v}(u^\mathrm{b})
  \quad \text{with} \quad
  u^\mathrm{b} = \frac{1}{a} g(\bm{x})
  \text{.}
\end{equation}

An important consideration is that boundary conditions are generally nonlinear.
Even for linear PDEs, such as the Poisson equation, a simple inhomogeneous
Dirichlet boundary condition $u^\mathrm{b} \neq 0$ introduces a nonlinearity in
the DG operator because $\dgA(0) \neq 0$. Therefore, we always linearize
boundary conditions. For nonlinear equations the boundary conditions linearize
along with the DG operator and require no special attention (see
\cref{sec:linearization}). However, for linear equations the inhomogeneity in
the boundary conditions is the only nonlinearity present in the DG operator, so
we skip the full linearization procedure. Instead, we contribute the
inhomogeneous boundary conditions~$\dgA(0)$ to the fixed sources, leaving only
the linearized boundary conditions in the DG operator,
\begin{equation}
  \frac{\delta\dgA}{\delta u} \grd{u} = \grd{\dgb} - \dgA(0)
  \text{,}
\end{equation}
where $\frac{\delta\dgA}{\delta u}$ is just $\dgA$ with linearized boundary
conditions. Note that this strategy is equivalent to the full linearization
procedure described in \cref{sec:linearization} at $\grd{u}=0$. In
practice, evaluating $\dgA(0)$ simplifies significantly for linear equations
because only the lifted external boundary corrections contribute to it.

\subsection{Evaluating the mass, stiffness, and lifting matrices}\label{sec:linop_impl}

The mass matrix~\eqref{eq:massmat}, the stiffness matrix~\eqref{eq:stiffmat}, and
the lifting operator~\eqref{eq:liftop} are integrals that must be evaluated on
every element of the computational domain. We evaluate these integrals on the
same grid on which we expand the dynamic fields, which amounts to a
Gauss-Lobatto quadrature of an order set by the number of collocation points in
the element. This strategy is commonly known as
\emph{mass lumping}.\footnote{This is the approach taken
in~\ccite{Teukolsky2016-ja}. See Eq.~(3.7) in~\ccite{Teukolsky2016-ja} for
details on the mass-lumped mass matrix on $d$-dimensional hexahedral elements.
Note that~\ccite{Teukolsky2016-ja} absorbs the metric determinant in the dynamic
variables.} Employing mass lumping and our choice of nodal
basis~\eqref{eq:lagrprod}, the mass matrix~\eqref{eq:massmat} evaluates to
\begin{equation}\label{eq:massmat_eval}
  \dgM_{pq} \approx \delta_{pq} \left.\sqrt{g}\right|_p \left.\jac\right|_p \prod_{i=1}^d w_{p_i}
  \text{.}
\end{equation}
Here the coefficients $w_{p_i}$ denote the Legendre-Gauss-Lobatto quadrature
weights associated with the collocation points $\xi_{p_i}$, and the geometric
quantities $\sqrt{g}$ and $\jac$ are evaluated directly on the collocation
points.\footnote{See, e.g., Algorithm~25 in \ccite{Kopriva} for details on
computing Legendre-Gauss-Lobatto quadrature weights, and Algorithm~23 for
Legendre-Gauss quadrature weights.} Recall from \cref{sec:domain} that the index
$p$ enumerates grid points in the regular $d$-dimensional grid and that $p_i$
denotes the grid point's index along the $i$th dimension. The diagonal
mass-lumping approximation~\eqref{eq:massmat_eval} has the advantage that it is
computationally efficient to apply, invert and store since it amounts to a
pointwise multiplication over the computational grid. Note that
\cref{eq:massmat_eval} is exact on a rectilinear LG grid with a flat background
geometry, and can be made exact on rectilinear LGL grids by including a
correction term without increasing the computational cost for applying or
inverting it \cite{Teukolsky2015-ov}. The quadrature weights~$w_{p_i}$ can be
cached and reused by all elements with the same number of collocation points in
a dimension.

The strong stiffness matrix~\eqref{eq:stiffmat} evaluates to
\begin{subequations}\label{eq:stiffmat_eval}
\begin{align}
  \dgMD_{i,pq} &\approx \dgM_{pr} \dgD_{i,rq}
  \text{,}
  \intertext{with}
  \dgD_{i,rq} &= \sum_{j=1}^d \left.\invjac^j_i\right|_r \lagr^\prime_{q_j}\!(\xi_{r_j})
  \prod_{\substack{k=1 \\ k\neq j}}^d\delta_{q_k r_k}
  \text{.}
\end{align}
\end{subequations}
Here $\ell^\prime_{q_j}\!(\xi_{r_j})$ are the one-dimensional \enquote{logical}
differentiation matrices obtained by differentiating the Lagrange
polynomials~\eqref{eq:lagr} and evaluating them at the collocation
points.\footnote{See, e.g., \ccite{HesthavenWarburton} and Algorithm~37
in~\ccite{Kopriva} for details on computing the one-dimensional logical
differentiation matrix $\lagr^\prime_{q_j}\!(\xi_{r_j})$ from properties of
Legendre polynomials.} The stiffness matrix is essentially a \enquote{massive}
differentiation operator that decomposes into one-dimensional differentiation
matrices due to the product structure of the basis functions~\eqref{eq:lagrprod}
on our hexahedral meshes. The one-dimensional logical differentiation matrices
can be cached and reused by all elements with the same number of collocation
points in a dimension, keeping the memory cost associated with the stiffness
operator to a minimum. The weak stiffness matrix can be computed analogously
from the transpose of the logical differentiation matrices.

The lifting operator~\eqref{eq:liftop} evaluates to
\begin{subequations}\label{eq:liftop_eval}
\begin{align}
  \dgML_{pq} &\approx \dgM_{pr} \dgL_{rq}
  \text{,}
  \intertext{with}
  \label{eq:liftop_eval_mat}
  \dgL_{rq} &= \delta_{rq} \sum_{i=1}^d (\delta_{{q_i} 1} + \delta_{{q_i}N_{k,i}})
  \frac{1}{w_{q_i}}
  \left.\sqrt{g^{jk}\invjac^i_j\invjac^i_k}\right|_{q}
  \mkern-8mu\text{.}
\end{align}
\end{subequations}
It is diagonal and has a contribution from every face of the element. Note that
each face only contributes to the LGL grid points on the respective face. On LG
grids additional interpolation matrices from the face into the volume appear in
this expression. Also note that the root in \cref{eq:liftop_eval_mat} is simply
the magnitude of the unnormalized face normal
$\tilde{n}_j$~\cite{Teukolsky2016-ja}.

Recall that the objective of these matrices is to evaluate the compact DG
operator \eqref{eq:dg_op_compact} along with \cref{eq:aux_matform} on every
element of the computational domain. In practice, neither matrix must be
assembled explicitly and stored on the elements: the mass
matrix~\eqref{eq:massmat_eval} reduces to a multiplication over the
computational grid, the stiffness matrix~\eqref{eq:stiffmat_eval} involves
logical one-dimensional differentiation matrices that are shared between
elements, and the lifting operation~\eqref{eq:liftop_eval} reduces to a
multiplication over the grid points on the element face. Since both the
stiffness and the lifting operation decompose into a mass matrix and a
\enquote{massless} operation, the same set of operations can be used to evaluate both
\cref{eq:aux_matform,eq:dg_op_compact}, and the mass matrix factors
out of the DG operator entirely. Nevertheless, we will see in \cref{sec:symm}
that it is advantageous to keep the mass matrix in the
operator~\eqref{eq:dg_op_compact}.

\subsection{A note on dealiasing}\label{sec:dealiasing}

The integral expressions discussed in \cref{sec:linop_impl} involve geometric
quantities that are typically known analytically, namely the Jacobian and the
background metric. Limiting the resolution of the integrals to the quadrature
order of the elements can make the scheme susceptible to geometric aliasing
because these quantities are resolved with limited precision, potentially
reducing the accuracy of the scheme on curved meshes or on a curved manifold. To
combat geometric aliasing we can, in principle, precompute the matrices on
every element at high accuracy, but at a significant memory cost. Precomputing
the matrices is possible in elliptic problems because the geometric quantities
remain constant. This is different to time-evolution systems that often involve
time-dependent Jacobians (\enquote{moving meshes}). Alternatively, a number of
dealiasing techniques are available to combat geometric aliasing, and also to
combat aliasing arising from evaluating other background quantities on the
collocation points, i.e., quantities in the PDEs that are independent of the
dynamic variables and known analytically~\cite{Mengaldo2015-qc}. For example,
\ccite{Vincent2019qpd} interpolates data from the primary LGL grid to an
auxiliary LG grid, on which the Jacobian is evaluated, to take advantage of the
higher-order quadrature. However, these dealiasing techniques can significantly
increase the computational cost for applying the DG operator. We have chosen to
employ the simple mass-lumping scheme detailed in \cref{sec:linop_impl} to
minimize the complexity, computational cost,  and memory consumption of the DG
operator. Detailed studies of cost-efficient dealiasing techniques are a
possible subject of future work.

\subsection{Mesh refinement}\label{sec:refinement}

The domain decomposition into elements, each with their own set of basis
functions, allows for two avenues to control the resolution: we can split the
domain into more and smaller elements ($h$~refinement) or increase the number of
basis functions within an element ($p$~refinement). We can perform $h$ and
$p$~refinement in each dimension independently.

\begin{figure}
  \centering
  \tikzinput{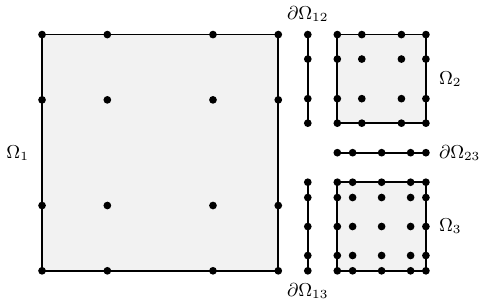}
  \caption{
    \label{fig:mortars}
    A representation of nonconforming element boundaries with mortars in two
    dimensions. The left element~$\Omega_1$ faces two elements $\Omega_2$ and
    $\Omega_3$ towards the right. Since both $\Omega_1$ and $\Omega_2$ have four
    vertical grid points their shared mortar~$\partial\Omega_{12}$ also has four
    grid points, but covers only a logical half of $\Omega_1$ ($h$~nonconforming).
    The element~$\Omega_3$ has five vertical grid points, so the
    mortar~$\partial\Omega_{13}$ has $\max(4, 5)=5$ grid points and also covers
    only a logical half of $\Omega_1$ ($hp$~nonconforming). Elements $\Omega_2$
    and~$\Omega_3$ are $h$~conforming but differ in their number of horizontal
    grid points, so their shared mortar~$\partial\Omega_{23}$ has $\max(4, 5)=5$
    grid points ($p$~nonconforming). Note that the empty space between the
    elements in this visualization is not part of the computational domain.}
\end{figure}

Both $h$~refinement and $p$~refinement can lead to nonconforming boundaries between
elements, meaning that grid points on the two sides of the boundary do not
coincide. Since we need to work with data from both sides of an element boundary
when considering numerical fluxes (see \cref{sec:numflux}) we place
\emph{mortars} between elements. A mortar is a $(d-1)$-dimensional mesh that has
sufficient resolution to exactly represent discretized fields from both adjacent
element faces. Specifically, a mortar $\partial\Omega_{k\bar{k}}$ between the
elements $\Omega_k$ and $\Omega_{\bar{k}}$ that share a boundary orthogonal to
dimension $j$ has $\max(N_{k,i},N_{\bar{k},i})$ grid points in dimension~$i\neq
j$. We limit the $h$~refinement of our computational domains such that an element
shares its boundary with at most two neighbors per dimension in every direction
(\enquote{two-to-one balance}). This means a mortar covers either the full element face
or a logical half of it in every dimension. \Cref{fig:mortars} illustrates an
$hp$-refined scenario with nonconforming element boundaries.

To project field values from an element face to a mortar we employ the
$(d-1)$-dimensional
\emph{prolongation operator}
\begin{equation}\label{eq:prol}
  \prol_{\tilde{p} p} = \prod_{i=1}^{d-1}\lagr_{p_i}(\tilde{\xi}_{\tilde{p}_i})
  \text{,}
\end{equation}
where $p$ enumerates grid points on the coarser (element face) mesh, $\tilde{p}$
enumerates grid points on the finer (mortar) mesh, and
$\tilde{\xi}_{\tilde{p}_i}$ are the coarse-mesh logical coordinates of the
fine-mesh collocation points. For mortars that cover the full element face in
dimension $i$ the coarse-mesh logical coordinates are just the fine-mesh
collocation points, $\tilde{\xi}_{\tilde{p}_i}=\xi_{\tilde{p}_i}$. For mortars
that cover the lower or upper logical half of the element face in dimension $i$
they are $\tilde{\xi}_{\tilde{p}_i}=(\xi_{\tilde{p}_i}-1)/2$ or
$\tilde{\xi}_{\tilde{p}_i}=(\xi_{\tilde{p}_i}+1)/2$, respectively. Note that the
prolongation operator~\eqref{eq:prol} is just a Lagrange interpolation from the
coarser (element face) mesh to the finer (mortar) mesh. The interpolation
retains the accuracy of the polynomial approximation because the mortar has
sufficient resolution. The prolongation operator is also an $L_2$~projection (or
\emph{Galerkin} projection) because it minimizes the $L_2$~norm
$\int_{\partial\Omega_{k\bar{k}}} (u^{(k)} - u^{(\tilde{k})})^2 \sqrt{g}
\dd{^{d-1}x}$, where $u^{(\tilde{k})}$ denotes the prolongated field values on
the finer (mortar) mesh.

To project field values from a mortar back to an element face we employ an
adjoint $\restr$ of the prolongation operator such that $\restr\prol=\ident$. We
also refer to this operation as a \emph{restriction} because it truncates higher
modes from the mortar down to the resolution of the element face. Specifically,
we employ the mass-conservative adjoint
\begin{align}
  \int_{\partial\Omega_{k\tilde{k}}} \mkern-10mu &\restr(u^{(\tilde{k})}) \, u^{(k)} \sqrt{g} \dd{^{d-1}x} \nonumber \\
  &= \int_{\partial\Omega_{k\tilde{k}}} \mkern-10mu u^{(\tilde{k})} \prol(u^{(k)}) \sqrt{g} \dd{^{d-1}x}
  \quad \forall u^{(k)}, u^{(\tilde{k})}
  \text{.}
\end{align}
In matrix notation the \emph{restriction operator} reduces to
\begin{equation}
  \restr = \dgM^{-1} \prol^T \tilde{\dgM}
  \text{,}
\end{equation}
where $\dgM^{-1}$ is the inverse mass matrix on the coarser (element face) mesh,
$\tilde{\dgM}$ is the mass matrix on the finer (mortar) mesh, and $\prol^T$ is
the transpose of the prolongation operator~\eqref{eq:prol}.

Note that the $d$-dimensional restriction and prolongation operators can serve
not only to project field values to and from mortars, but also to project field
values to and from elements that cover the computational domain at different
$h$- and $p$-refinement levels. We make no use of projections across refinement levels
in this article but will do so in upcoming work for the purpose of adaptive
mesh-refinement strategies and for multigrid solvers.\footnote{See also
Sections~3.2 and~3.3 in \ccite{Fortunato2019} for details on the restriction
and prolongation operators in the context of multigrid solvers.}

\subsection{A note on symmetry}\label{sec:symm}

For practical applications it is often advantageous to work with a symmetric
operator. For example, some iterative linear solvers such as conjugate gradients
take advantage of the symmetry to invert the operator more efficiently. One can
also often show stronger convergence bounds for iterative linear solvers
applicable to nonsymmetric matrices, such as GMRES, if the matrix is
symmetric~\cite{Saad1986-ix}.

The compact strong-weak scheme presented in \cref{eq:dg_op_compact_weak} with
the generalized SIP numerical flux~\eqref{eq:numflux} is symmetric unless the
elliptic equations break the symmetry, e.g., with an asymmetric coupling between
equations. Note that a curved manifold will typically break the symmetry because
it involves first derivatives in Christoffel-symbol contributions to the primal
sources [see, e.g., \cref{eq:poisson_curved_fluxform}]. It is straightforward to
see how the strong-weak scheme can make the DG operator symmetric if the
elliptic equations allow it: the strong-weak operator involves a symmetric
stiffness term of the schematic form $(\dgM \dgD)^T \dgD = \dgD^T \dgM \dgD$,
whereas the strong scheme has a nonsymmetric expression of the form $\dgM \dgD
\dgD$ instead. Note that the \enquote{massless} variant of the strong-weak scheme,
schematically~$\dgM^{-1} \dgD^T \dgM \dgD$, is not generally symmetric, and
neither is the \enquote{massless} strong scheme~$\dgD \dgD$.

\subsection{Linearizing the operator}\label{sec:linearization}

To solve nonlinear elliptic equations $\dgA(\grd{u}) = \grd{\dgb}$
we typically employ a correction scheme, repeatedly solving the linearized
equations for a correction quantity~$\Delta \grd{u}$. For example, a
simple Newton-Raphson correction scheme solves the linearized problem
$\frac{\delta\dgA}{\delta u}(\grd{u}) \, \Delta \grd{u} =
\grd{\dgb} - \dgA(\grd{u})$ at fixed $\grd{u}$ and then
iteratively corrects $\grd{u} \rightarrow \grd{u} + \Delta
\grd{u}$. Since the fluxes $\dgFi{\alpha}$ are already linear for all
elliptic systems we consider, the linearization $\frac{\delta\dgA}{\delta
u}(\grd{u})$ involves only linearizing the sources $\dgS_\alpha$ and the
boundary conditions.

\subsection{Variations of the scheme}\label{sec:variations}

We have made a number of choices to formulate the DG discretization in the
preceding sections. This section summarizes some of the choices and presents
possible variations to explore in future work.

\begin{description}
  \item[Massive vs massless scheme] We can eliminate the mass matrix
  in \cref{eq:dg_op_compact} to obtain a \enquote{massless} DG operator. However, we
  have found evidence that iterative linear solvers converge faster when solving
  the \enquote{massive} DG operator. We attribute this behaviour to the symmetry
  considerations discussed in \cref{sec:symm}.
  \item[Mass-lumping] We diagonally approximate the mass matrix to reduce the
  computational cost to apply, invert and store it, and to simplify the scheme
  (see \cref{sec:linop_impl}). Dealiasing techniques can potentially
  increase the accuracy of the scheme on curved meshes as discussed
  in \cref{sec:dealiasing}.
  \item[LGL vs LG mesh] We chose to discretize the DG operator on LGL meshes
  to take advantage of the collocation points on element boundaries, which
  simplify computations of boundary corrections. Switching to LG meshes can have
  the advantage that quadratures are one degree more precise, making the
  mass-lumping exact on rectilinear grids (see \cref{sec:linop_impl}).
  \item[Numerical flux] The generalized internal-penalty numerical flux
  presented in \cref{sec:numflux} has proven a viable choice for a wide range of
  problems so far. However, the ability to switch out the numerical flux is a
  notable strength of DG methods, and augmenting the numerical flux in the
  elliptic DG scheme may improve its convergence properties or accuracy. In
  particular, the choice of penalty, \cref{eq:penalty}, on curved meshes remains
  a subject of further study.
  \item[Strong vs weak formulation] We have chosen the strong
  formulation~\eqref{eq:dg_op_compact_strong} over the strong-weak
  formulation~\eqref{eq:dg_op_compact_weak} because it is slightly simpler and
  we have, so far, found no evidence that the strong-weak formulation converges
  faster than the strong formulation, despite the symmetry considerations
  discussed in \cref{sec:symm}. However, the strong-weak formulation can be
  of interest if a symmetric DG operator is necessary, e.g., to take advantage of
  specialized iterative solvers.
  \item[Flux vs primal formulation] We have eliminated auxiliary degrees of
  freedom in the DG operator with a Schur-complement strategy. An alternative
  strategy is to derive a \enquote{primal formulation} of the DG operator
  (see \cref{sec:compact}). We have found the flux formulation easier to
  implement due to its similarity to hyperbolic DG schemes. Furthermore,
  \ccite{Fortunato2019} suggests that the flux formulation can be advantageous in
  conjunction with a multigrid solver.
\end{description}

\section{Test problems}\label{sec:tests}

The following numerical tests confirm the DG scheme presented in this article
can solve a variety of elliptic problems. The test problems involve linear and
nonlinear systems of PDEs with nonlinear boundary conditions on curved
manifolds, discretized on $hp$-refined domains with curved meshes and
nonconforming element boundaries.

For test problems that have an analytic solution we quantify the accuracy of the
numerical solutions by computing an $L_2$~error over all primal variables,
\begin{equation}\label{eq:errnorm}
  \lVert u-u_\mathrm{analytic} \rVert \defeq \left(\frac{
    \sum_{A,k} \int_{\Omega_k} (u_A - u_{A,\mathrm{analytic}})^2 \dd{V}}{
      \sum_k \int_{\Omega_k} \dd{V}}
    \right)^{1/2}
    \mkern-22mu \text{,}
\end{equation}
where the integrals are evaluated with Gauss-Lobatto quadrature on the elements
of the computational domain.

To assess the DG operator is functional for our test problems we study the
convergence of the discretization error~\eqref{eq:errnorm} under uniform
$hp$~refinement of the computational domain (see \cref{sec:refinement}). We
compute the $h$-convergence order under pure uniform $h$~refinement
\begin{equation}
  \tau_h \defeq \frac{\Delta_h\ln(\lVert u - u_\mathrm{analytic}\rVert)}{\Delta_h\ln(h)}
  \text{,}
\end{equation}
where $\Delta_h$ denotes the difference between successive $h$-refinement levels
and $h$ is the size of an element. Since we always split elements in half along
all logical axes we use $\Delta_h\ln(h)=\ln(2)$. We also compute the exponential
convergence scale under pure uniform $p$~refinement
\begin{equation}
  \tau_p \defeq \Delta_p\log_{10}(\lVert u - u_\mathrm{analytic}\rVert)
  \text{,}
\end{equation}
where $\Delta_p$ denotes the difference between successive $p$-refinement levels.

\subsection{A Poisson solution}\label{sec:test_poisson}

\begin{figure}
  \centering
  \includegraphics[width=0.6\columnwidth]{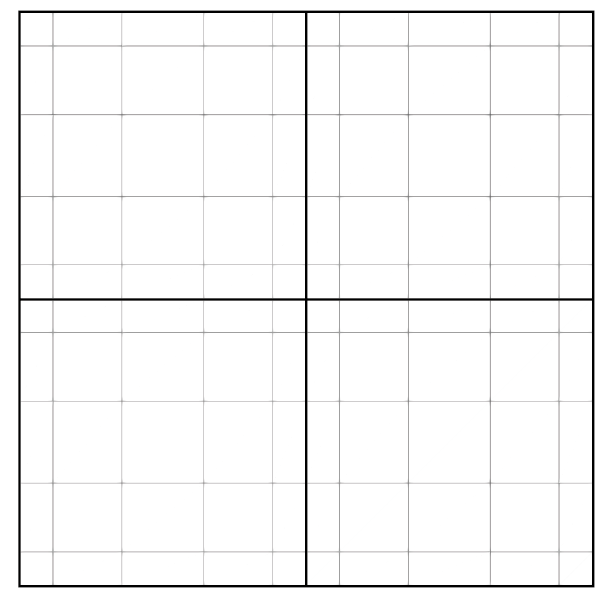}
  \caption{
    \label{fig:poisson_domain}
    The two-dimensional rectilinear domain used in the Poisson problem
    (\cref{sec:test_poisson}). Black lines illustrate element boundaries and
    gray lines represent the LGL grid within each element. This domain is
    isotropically $h$~refined once, i.e., split once in both dimensions, resulting
    in four elements. Each element has six grid points per dimension, so fields
    are represented as polynomials of degree five. This is the domain that
    \cref{fig:noncompact_operator,fig:compact_operator} are based on,
    and that is circled in \cref{fig:poisson_convergence}.}
\end{figure}

\begin{figure*}
  \centering
  \includegraphics[width=1.5\columnwidth]{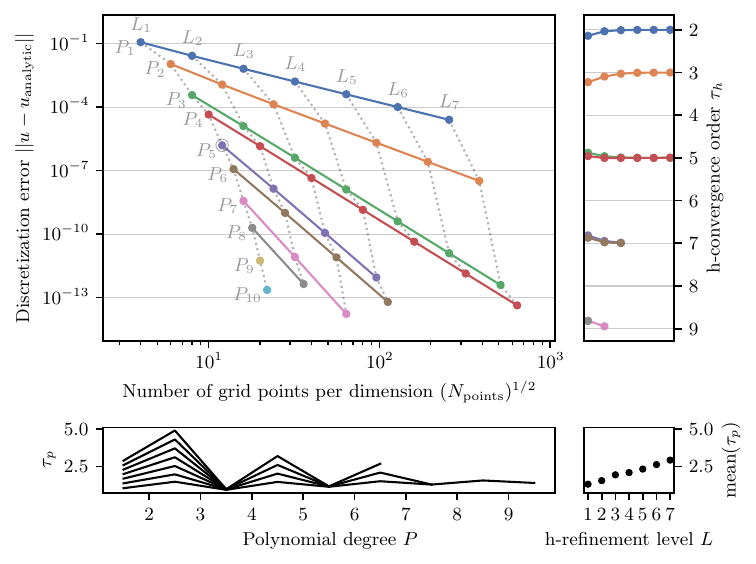}
  \caption{
    \label{fig:poisson_convergence}
    Convergence of the two-dimensional Poisson problem detailed in
    \cref{sec:test_poisson} with uniform $hp$~refinement. Solid lines connect
    numerical solutions where the domain is split into an increasing number of
    elements (isotropic $h$~refinement), and dotted lines connect numerical
    solutions with increasing polynomial order (isotropic $p$~refinement). The DG
    scheme recovers optimal $\mathcal{O}(h^{P+1})$~convergence with odd-order
    superconvergence under $h$~refinement (right panel) and exponential
    convergence under $p$~refinement (bottom panels). For reference, the circled
    configuration is pictured in \cref{fig:poisson_domain}.}
\end{figure*}

With this first test problem we establish a simple baseline that the following
tests build upon. It is reduced to the absolute essentials to illustrate the
basic concepts of the scheme. We solve a flat-space Poisson
equation~\eqref{eq:poisson} in two dimensions for the analytic solution
\begin{equation}\label{eq:prodofsin}
  u_\text{analytic}(\bm{x})=\sin{(\pi x)}\sin{(\pi y)}
\end{equation}
on a rectilinear domain $\Omega=[0,1]^2$. The domain is illustrated in
\cref{fig:poisson_domain}. To obtain the solution~\eqref{eq:prodofsin}
numerically we choose the fixed source~$\dgf(\bm{x})=2\pi^2\sin{(\pi
x)}\sin{(\pi y)}$, select homogeneous Dirichlet boundary conditions
${u^\mathrm{b}=0}$, and solve the strong compact DG-discretized
problem~\eqref{eq:dg_op_compact_strong} with~${C=1}$. This essentially means we
invert the matrix depicted in \cref{fig:compact_operator} and apply it to the
discretization of the fixed source~$\dgf(\bm{x})$. Instead of inverting the
matrix directly we employ the iterative elliptic solver of the \spectre{}
code~\cite{spectre} presented in~\ccite{ellsolver}. However, note that the
technology we use to solve the DG-discretized problem is not relevant for the
purpose of this article, since the matrix equation has a unique solution.
Assuming the matrix equation is solved to sufficient precision,
\cref{eq:errnorm} quantifies the discretization error of the DG scheme.

We solve the problem on a series of uniformly and isotropically refined domains
and present the convergence of the discretization error in
\cref{fig:poisson_convergence}. Under $h$~refinement the scheme recovers
optimal $\mathcal{O}(h^{P+1})$~convergence, where $P$ denotes the polynomial
degree of the elements. It also recovers the odd-order superconvergence feature
expected for the antisymmetric problem~\eqref{eq:prodofsin}.\footnote{See also
Fig.~7.9 in~\ccite{HesthavenWarburton}.} Under $p$~refinement the scheme recovers
exponential convergence. The exponential convergence scale~$\tau_p$ is modulated
by the superconvergence feature and its mean increases linearly with the
$h$-refinement level.

\subsection{Thermal noise in a cylindrical mirror}\label{sec:test_elasticity}

\begin{table}
\begin{ruledtabular}
\begin{tabular}{llr}
  Beam width & $r_0$ & \SI{177}{\micro\metre} \\
  Outer radius & $R$ & \SI{600}{\micro\metre} \\
  \hline
  Poisson ratio & $\nu$ & 0.17 \\
  Young's modulus & $E$ & \SI{72}{GPa} \\
\end{tabular}
\caption{
  \label{tab:elasticityparams}
  Parameters used in the thermal-noise problem (\cref{sec:test_elasticity}).
  The beam width and the material properties correspond to Table 1
  in~\ccite{Lovelace2017xyf}. These material properties characterize a
  fused-silica mirror, which is a material used in the LIGO gravitational-wave
  detectors.}
\end{ruledtabular}
\end{table}

\begin{figure}
  \centering
  \includegraphics[width=0.9\columnwidth]{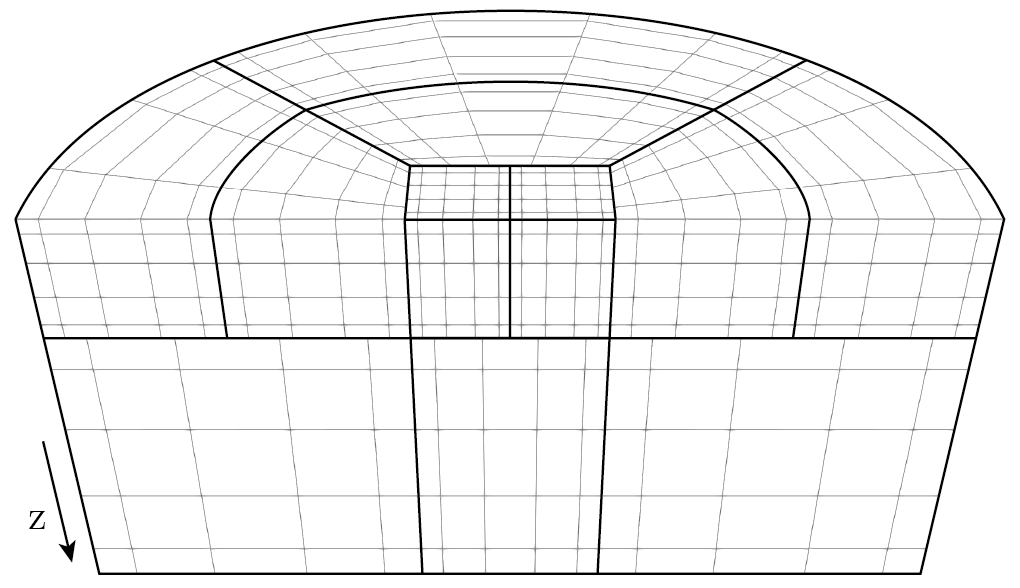}
  \caption{
    \label{fig:elasticity_domain}
    A cut through the cylindrical domain used in the elasticity problem
    (\cref{sec:test_elasticity}). The domain consists of four wedges enveloping
    a cuboid, and two vertical layers. The layers are partitioned vertically at
    $z=r_0$ and the cuboid lies radially within $r=r_0$. In the top layer, the
    wedges are $h$~refined radially once and the cuboid is $h$~refined in the
    $x$~and $y$~directions once, resulting in 12 elements in the top layer and 5
    elements in the bottom layer. Elements in this example have six grid points
    per dimension, and the wedge-shaped elements have two additional grid points
    in their angular direction.}
\end{figure}

\begin{figure*}
  \centering
  \includegraphics[width=1.5\columnwidth]{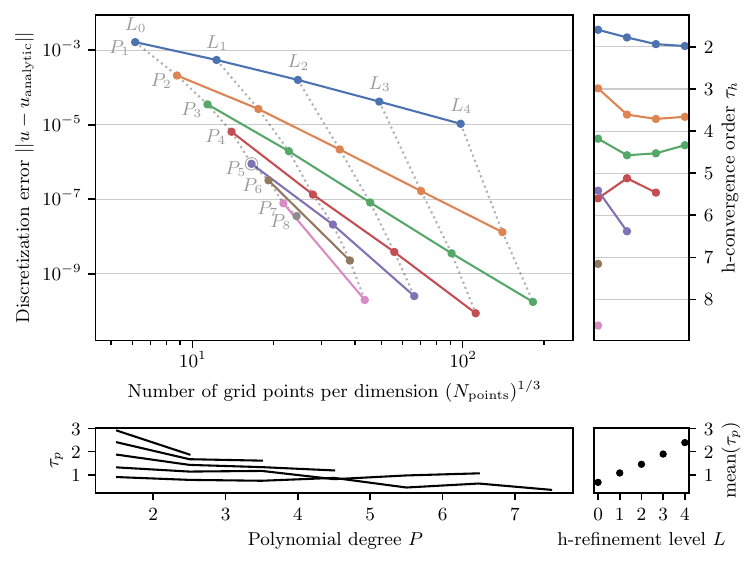}
  \caption{
    \label{fig:elasticity_convergence}
    Convergence of the three-dimensional elasticity problem detailed in
    \cref{sec:test_elasticity} under uniform $hp$~refinement. Plotted is the
    $L_2$~error~\eqref{eq:errnorm} over the three components of the displacement
    field~$\xi^i(\bm{x})$. The refinement is based on the domain pictured in
    \cref{fig:elasticity_domain} (circled configuration) with curved meshes and
    nonconforming element boundaries. The DG scheme recovers
    optimal $\mathcal{O}(h^{P+1})$~convergence under $h$~refinement (right panel)
    and exponential convergence under $p$~refinement (bottom panels).}
\end{figure*}

In this second test problem we solve the equations of linear
elasticity~\eqref{eq:elasticity} on a curved mesh with nonconforming element
boundaries. The test problem represents a cylindrical mirror that is deformed by
pressure from a laser beam incident on one of the sides. This problem arises in
studies of Brownian thermal noise in interferometric gravitational-wave
detectors \cite{Levin:1997kv,Lovelace2017xyf}.\footnote{See also Section~11.9.2
in \ccite{BlandfordThorne} for an introduction to the thermal noise problem.}
Here we consider an analytic solution to this problem that applies in the limit
of an isotropic and homogeneous mirror material with constitutive relation
\begin{equation}\label{eq:constrel_isotropic_homogeneous}
  \constrel^{ijkl} = \lambda \, \delta^{ij} \delta^{kl} + \mu \left(\delta^{ik}\delta^{jl} + \delta^{il}\delta^{jk}\right)
  \text{,}
\end{equation}
characterized by the Lam\'{e} parameter~$\lambda$ and the shear modulus~$\mu$, or
equivalently by the Poisson ratio~$\nu=\frac{\lambda}{2(\lambda+\mu)}$, Young's
modulus~$E=\frac{\mu(3\lambda+2\mu)}{\lambda+\mu}$, or the bulk modulus
$K=\lambda + \frac{2}{3}\mu$. We assume the material fills the infinite
half-space $z \geq 0$, choose a vanishing force density~$f^j(x)=0$, and a
Gaussian profile of the laser beam incident at $z=0$,
\begin{equation}\label{eq:laserbeam}
  n_i \stress^{ij} = n^j\frac{1}{\pi r_0^2} e^{-r^2 / r_0^2}
  \text{.}
\end{equation}
Here $\stress^{ij}=-\constrel^{ijkl} S_{kl}$ is the \emph{stress}, $n_i$ is the
unit normal pointing away from the mirror, i.e., in negative $z$~direction,
$r=\sqrt{x^2 + y^2}$ is the radial coordinate distance from the axis of
symmetry, and $r_0$ is the beam width. Under these assumptions the displacement
field~$\xi^i(\bm{x})$ has the analytic
solution~\cite{Liu2000-qz,BlandfordThorne,Lovelace2007tn}
\begin{subequations}\label{eq:halfspacemirror}
\begin{align}
  \displ^r &= \frac{1}{2\mu} \int_0^{\infty} \! \dd{k} \besselJ_1(kr)e^{-kz}
    \left(1 - \frac{\lambda + 2\mu}{\lambda + \mu} + kz \right) \tilde{p}(k)
    \text{.} \\
    \displ^z &= \frac{1}{2\mu} \int_0^{\infty} \! \dd{k} \besselJ_0(kr)e^{-kz}
      \left(1 + \frac{\mu}{\lambda + \mu} + kz \right) \tilde{p}(k)
\end{align}
\end{subequations}
and $\displ^\phi=0$ in cylindrical coordinates $\{r,\phi,z\}$. Here
$\besselJ_0$ and $\besselJ_1$ are Bessel functions of the first kind, and
$\tilde{p}(k) = \frac{1}{2\pi} e^{-(k r_0 / 2)^2}$ is the Hankel transform of
the laser-beam profile. We evaluate these integrals numerically at every
collocation point in the computational domain to determine the analytic
solution.

To obtain numerical solutions to the thermal noise problem we DG discretize the
equations of linear elasticity~\eqref{eq:elasticity_fluxform} on a cylindrical
domain with height and radius~$R$, employing the strong compact DG
operator~\eqref{eq:dg_op_compact_strong}. Since the
stress~$T^{ij}=-\tensor{\dgF}{_{\!\displ}^i^j}$ is the negative primal flux in
the elasticity equations~\eqref{eq:elasticity_fluxform} we
impose \cref{eq:laserbeam} as a Neumann-type boundary condition on the base of
the cylinder at $z=0$. We impose the analytic
solution~\eqref{eq:halfspacemirror} as Dirichlet-type boundary conditions on the
remaining external boundaries of the domain, i.e., on the base at $z=R$ and on
the mantle at $r=R$. These boundary conditions mean that we solve for a finite
cylindrical section of the infinite half-space analytic
solution~\eqref{eq:halfspacemirror}. We choose a penalty parameter of~${C=100}$
for this problem to eliminate variations in the discretization error arising
from curved-mesh contributions to the penalty~\eqref{eq:penalty} at high
resolutions. \Cref{tab:elasticityparams} summarizes the remaining parameters
we use in the numerical solutions.

\Cref{fig:elasticity_domain} illustrates the cylindrical domain. It is refined
more strongly toward the origin~$\bm{x}=0$ where the Gaussian laser beam applies
pressure. The refinement is both anisotropic and inhomogeneous, leading to
nonconforming element boundaries with different polynomial degrees on either
side of the boundary, multiple neighbors adjacent to an element face, or both.
Specifically, elements facing the top-layer cuboid or the interface between top
and bottom layer are matched two-to-one, and wedge-shaped elements have two
additional angular grid points. Therefore, the elements facing the cuboid are
both $p$~nonconforming and $h$~nonconforming in the top layer, and $p$~nonconforming
in the bottom layer. The elements facing the layer interface are
$h$~nonconforming.

\Cref{fig:elasticity_convergence} presents the convergence of the discretization
error under uniform $hp$~refinement. Specifically, we split every element in two
along all three dimensions to construct additional $h$-refinement levels, and
increment every polynomial degree by one to construct additional $p$-refinement
levels, retaining the nonconforming element boundaries. Note that the
wedge-shaped elements retain a higher polynomial degree of $P+2$ along their
angular direction throughout the refinement procedure, where $P$ is the
polynomial degree of all other elements and dimensions. The DG scheme recovers
optimal $\mathcal{O}(h^{P+1})$~convergence under $h$~refinement and exponential
convergence under $p$~refinement. Note that the exponential convergence
scale~$\tau_p$ depends on the domain geometry, the structure of the solution,
the placement of grid points and the refinement strategy. We have chosen to
refine the domain as uniformly as possible here to reliably measure convergence
properties of the DG scheme. Optimizing the distribution of elements and grid
points with adaptive mesh refinement (AMR) strategies to increase the rate of
convergence is the subject of ongoing work.

\subsection{A black hole in general relativity}\label{sec:kerrschild}

\begin{figure}
  \centering
  \includegraphics[width=0.9\columnwidth]{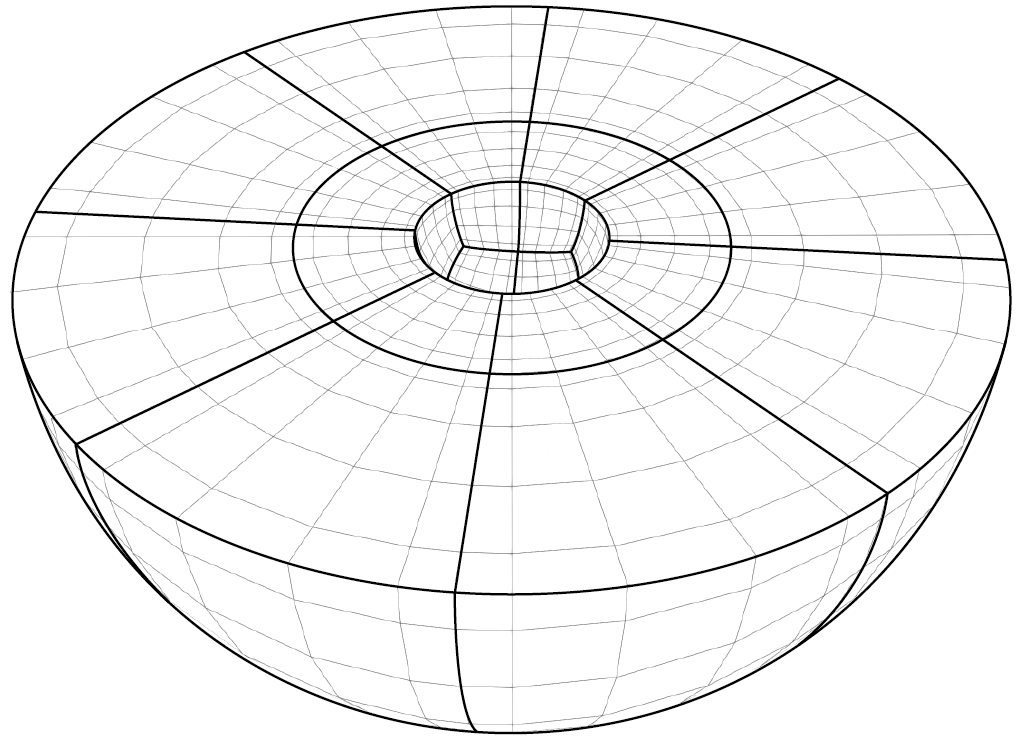}
  \caption{
  \label{fig:kerrschild_domain}
    A cut through the uniformly refined spherical-shell domain used in the
    black hole problem (\cref{sec:kerrschild}). The domain consists of six
    wedges with a logarithmic radial coordinate map enveloping an excised
    sphere. In this example each wedge is isotropically $h$~refined once, i.e.,
    split once in all three dimensions, resulting in a total of 48 elements.
    Note the elements are split in half along their logical axes, so the element
    size scales logarithmically in radial direction just like the distribution
    of grid points within the elements. Each element has six grid point per
    dimension, so fields are represented as polynomials of degree five.}
\end{figure}

\begin{figure*}
  \centering
  \includegraphics[width=1.5\columnwidth]{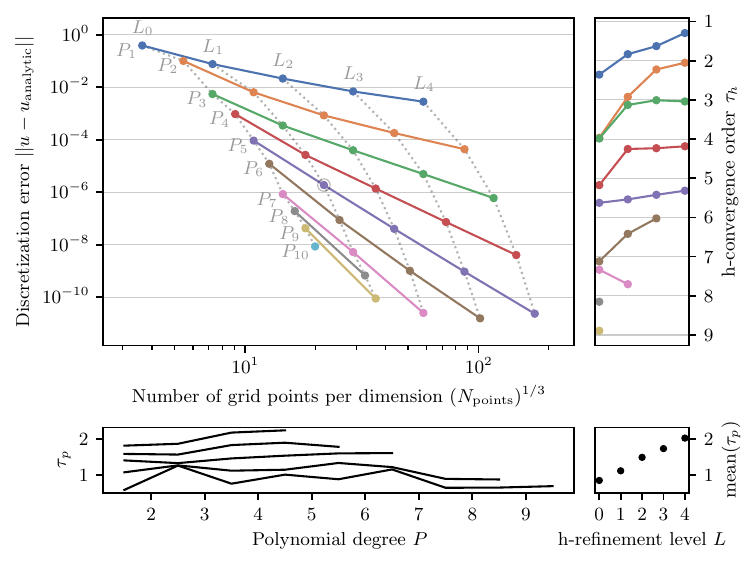}
  \caption{
  \label{fig:kerrschild_convergence}
    Convergence of the three-dimensional black-hole solution with uniform
    $hp$~refinement. Plotted is the $L_2$~error~\eqref{eq:errnorm} over all
    variables of the XCTS equations $\{\psi,\alpha\psi,\beta^i\}$. The circled
    configuration is pictured in \cref{fig:kerrschild_domain}. The DG scheme
    recovers $\mathcal{O}(h^P)$~convergence under $h$~refinement (right panel) and
    exponential convergence under $p$~refinement (bottom panels).}
\end{figure*}

Now we apply the DG scheme to solve the Einstein constraint equations of
general relativity in the XCTS formulations, which is a set of coupled,
nonlinear, elliptic PDEs on a curved manifold (see \cref{sec:xcts}). Solutions
to the XCTS equations describe admissible configurations of general-relativistic
spacetime and provide initial data for general-relativistic time evolutions.

In this test problem we solve the XCTS equations~\eqref{eq:xcts_fluxform} for a
Schwarzschild black hole in Kerr-Schild coordinates,
\begin{subequations}\label{eq:kerrschild}
\begin{align}
  \label{eq:kerrschild_psi}
  \psi &= 1
  \text{,} \\
  \label{eq:kerrschild_lapse}
  \alpha &= \left(1 + \frac{2M}{r}\right)^{-1/2}
  \text{,} \\
  \label{eq:kerrschild_shift}
  \beta^i &= \frac{2M}{r} \alpha^2 l^i
  \text{,}
  \intertext{with the background quantities}
  \bar{\gamma}_{ij} &= \delta_{ij} + \frac{2M}{r} l_i l_j
  \intertext{and}
  K &= \frac{2M\alpha^3}{r^2}\left(1 + \frac{3M}{r}\right)
  \text{,}
\end{align}
\end{subequations}
where $M$ is the mass parameter, $r=\sqrt{x^2+y^2+z^2}$ is the Euclidean
coordinate distance, and $l^i=l_i=x^i/r$.\footnote{See Table 2.1
in~\ccite{BaumgarteShapiro}.} The time-derivative quantities $\bar{u}_{ij}$ and
$\partial_t K$ in the XCTS equations~\eqref{eq:xcts_fluxform} vanish, as do the
matter sources $\rho$, $S$, and~$S^i$. Note that we have chosen a conformal
decomposition with $\psi=1$ here, but other choices of $\psi$ and
$\bar{\gamma}_{ij}$ that keep the spatial metric $\gamma_{ij} = \psi^4
\bar{\gamma}_{ij}$ invariant are equally admissible.

We solve the XCTS equations numerically for the conformal factor~$\psi$, the
product~$\alpha\psi$, and the shift~$\beta^i$. The conformal
metric~$\bar{\gamma}_{ij}$ and the trace of the extrinsic curvature $K$ are
background quantities that are chosen in advance and remain fixed throughout the
solve. They are a source of aliasing when evaluated on the computational grid
(see \cref{sec:dealiasing}). Importantly for this test problem the conformal
metric $\bar{\gamma}_{ij}$ is not flat, resulting in a problem formulated on a
curved manifold. For example, unit normal one forms in the DG scheme are
normalized with respect to the conformal metric $\bar{\gamma}_{ij}$ and the
metric determinant appears in the mass matrix and in the
$L_2$~error~\eqref{eq:errnorm}.

To solve the black hole problem numerically we employ the strong compact DG
scheme~\eqref{eq:dg_op_compact_strong} with~${C=1}$ to discretize the XCTS
equations~\eqref{eq:xcts_fluxform} on a three-dimensional spherical shell, as
illustrated in \cref{fig:kerrschild_domain}. The domain envelops an excised
sphere that represents the black hole, so it has an outer and an inner external
boundary that require boundary conditions. To obtain the Schwarzschild solution
in Kerr-Schild coordinates we impose
\crefrange{eq:kerrschild_psi}{eq:kerrschild_shift} as Dirichlet-type boundary
conditions at the outer boundary of the spherical shell at~$r=10M$. We place the
inner radius of the spherical shell at $r=2M$ and impose nonspinning
apparent-horizon boundary conditions at the inner boundary,
\begin{subequations}\label{eq:ah_bc}
\begin{align}
  n^k\partial_k\psi = &\frac{\psi^3}{8\alpha}n_i n_j\left((\bar{L}\beta)^{ij} - \bar{u}^{ij}\right) \nonumber \\
  &- \frac{\psi}{4}\bar{m}^{ij}\bar{\nabla}_in_j - \frac{1}{6}K\psi^3 \text{,} \\
  \beta^i = &-\frac{\alpha}{\psi^2}n^i
  \text{,}
\end{align}
\end{subequations}
where $\bar{m}^{ij}=\bar{\gamma}^{ij}-n^i n^j$. These boundary conditions are
not specific to the Schwarzschild solution but ensure the excision surface is an
\emph{apparent horizon}~\cite{Cook2004-yf}.\footnote{See e.g., Section~12.3.2
in~\ccite{BaumgarteShapiro} for an introduction to the apparent-horizon boundary
conditions.} Since the Schwarzschild solution in Kerr-Schild coordinates has an
apparent horizon at $r=2M$ we recover the solution~\eqref{eq:kerrschild} when we
place the inner radius of the spherical shell at that radius. The
apparent-horizon boundary conditions~\eqref{eq:ah_bc} do not constrain the
lapse~$\alpha$, so we impose \cref{eq:kerrschild_lapse} at the inner boundary.
The apparent-horizon boundary conditions are of Neumann-type for the variable
$\psi$, of Dirichlet-type for $\alpha\psi$ and $\beta^i$, and nonlinear.

Since the XCTS equations~\eqref{eq:xcts} and the apparent-horizon boundary
conditions~\eqref{eq:ah_bc} are nonlinear the initial guess for the iterative
nonlinear solver becomes relevant. We choose an initial guess close to the
analytic solution to ensure fast convergence of the iterative solver to the
numerical solution. Note that the initial guess and other details of the
iterative solve do not affect the discretization error of the numerical solution
once the solve has converged to sufficient precision.

We present the convergence of the discretization error under uniform
$hp$~refinement in \cref{fig:kerrschild_convergence}. The DG scheme for the
nonlinear black hole problem recovers $\mathcal{O}(h^P)$~convergence under
$h$~refinement, which is an order lower than that obtained for the two preceding linear
test problems. We find higher-order convergence for pure Dirichlet boundary
conditions for this problem, suggesting the apparent-horizon boundary
conditions~\eqref{eq:ah_bc} are responsible for the reduction of the convergence
order. For a Poisson problem with nonlinear boundary conditions the authors of
\ccite{Feistauer2019-cw} also find a loss of convergence under $h$~refinement.
Under $p$~refinement the scheme recovers exponential convergence and the mean
exponential convergence scale~$\tau_p$ increases linearly with the $h$-refinement
level.

\subsection{A black hole binary}\label{sec:bbh}

Finally, we solve the Einstein constraint equations in the XCTS formulation as
in \cref{sec:kerrschild}, but now we choose background quantities and boundary
conditions that represent two black holes in orbit. This binary black hole
problem is of significant relevance in numerical relativity to procure initial
data for simulations of merging black holes~\cite{BaumgarteShapiro,
Pfeiffer2004-oo, Lovelace2008-sw, Varma2018-fp}.

Following the formalism for \emph{superposed Kerr-Schild} initial data, e.g.,
laid out in \ccite{Lovelace2008-sw, Varma2018-fp}, we set the conformal metric
and the trace of the extrinsic curvature to the superpositions
\begin{subequations}\label{eq:bbh_bg}
\begin{align}
  \bar{\gamma}_{ij} &= \delta_{ij} + \sum_{n=1}^2 e^{-r_n^2 / w_n^2} \, (\gamma_{ij}^{(n)} - \delta_{ij})
\intertext{and}
  K &= \sum_{n=1}^2 e^{-r_n^2 / w_n^2} K^{(n)}
  \text{,}
\end{align}
\end{subequations}
where $\gamma_{ij}^{(n)}$ and $K^{(n)}$ are the conformal metric and
extrinsic-curvature trace of two isolated Schwarzschild black holes in
Kerr-Schild coordinates as given in Eqs.~(\ref{eq:kerrschild}). They have mass parameters~$M_n$
and are centered at coordinates~$\bm{C}_n$, with $r_n$ being the Euclidean coordinate
distance from either center. The superpositions are modulated by two Gaussians
with widths $w_n$. The time-derivative quantities $\bar{u}_{ij}$ and $\partial_t
K$ in the XCTS equations~\eqref{eq:xcts} vanish, as do the matter sources
$\rho$, $S$ and~$S^i$.

\begin{figure}
  \centering
  \includegraphics[width=0.9\columnwidth]{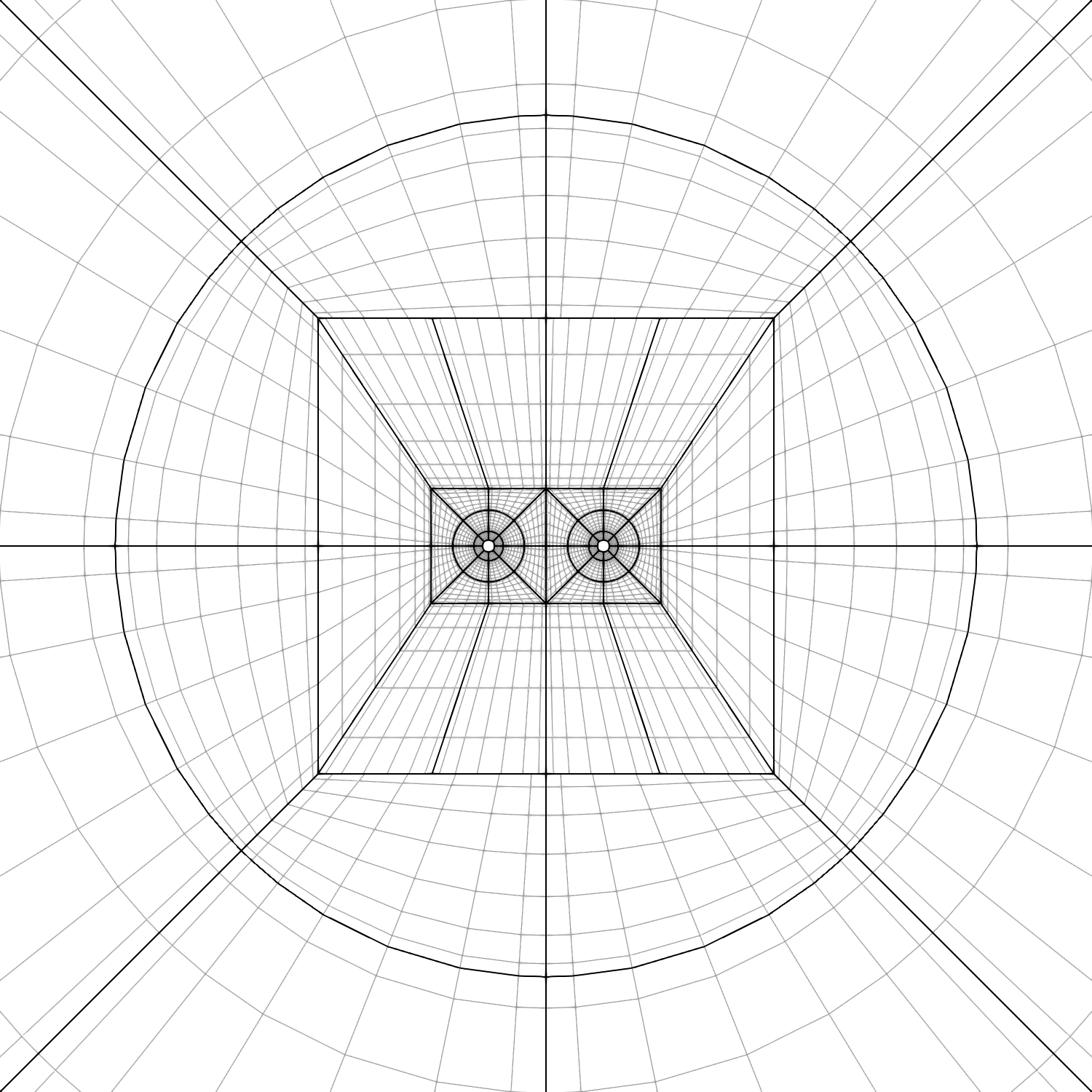}
  \caption{
    \label{fig:bbh_domain}
    A cut through the three-dimensional black-hole binary domain used in
    \cref{sec:bbh}. It involves two excised spheres centered at $\bm{C}_n$ along
    the $x$~axis and extends to a spherical outer surface at radius~$R$. The
    domain is $h$~refined such that spherical wedges have equal angular size, so
    the cube-to-sphere boundary is nonconforming. All elements in this picture
    have eight angular grid points, and $\{7, 8, 8, 9, 11, 11\}$ radial grid
    points in the layers ordered from outermost to innermost.}
\end{figure}

To handle orbital motion we split the shift in a \emph{background} and an
\emph{excess} contribution~\cite{Pfeiffer2003-Thesis},
\begin{equation}\label{eq:shift_split}
  \beta^i = \beta^i_\mathrm{background} + \beta^i_\mathrm{excess}
  \text{,}
\end{equation}
and choose the background shift
\begin{equation}\label{eq:bbh_bg_shift}
  \beta^i_\mathrm{background} = (\bm{\Omega}_0 \times \bm{x})^i
  \text{,}
\end{equation}
where $\bm{\Omega}_0$ is the orbital angular velocity. We insert
\cref{eq:shift_split} in the XCTS equations~\eqref{eq:xcts} and henceforth solve
them for $\beta^i_\mathrm{excess}$, instead of $\beta^i$.

We solve the XCTS equations on the domain depicted in \cref{fig:bbh_domain}. It
has two excised spheres with radius $2 M_n$ that are centered at $\bm{C}_n$, and
correspond to the two black holes, and an outer spherical boundary at finite
radius~$R$. We impose boundary conditions on these three boundaries as follows.
At the outer spherical boundary we impose asymptotic flatness,
\begin{equation}\label{eq:flatness_bc}
  \psi = 1 \text{,} \quad \alpha \psi = 1 \text{,} \quad \beta_\mathrm{excess}^i = 0
  \text{.}
\end{equation}
Since the outer boundary is at a finite radius, the solution will only be
approximately asymptotically flat.
On the two excision boundaries we impose nonspinning apparent-horizon boundary
conditions, \cref{eq:ah_bc}. For the lapse we choose to impose the isolated
solution~\eqref{eq:kerrschild_lapse} as Dirichlet conditions at both excision
surfaces. Note that this choice differs slightly from \ccite{Varma2018-fp},
where the \emph{superposed} isolated solutions are imposed on the lapse at both
excision surfaces. 

\begin{figure}
  \centering
  \includegraphics[width=\columnwidth]{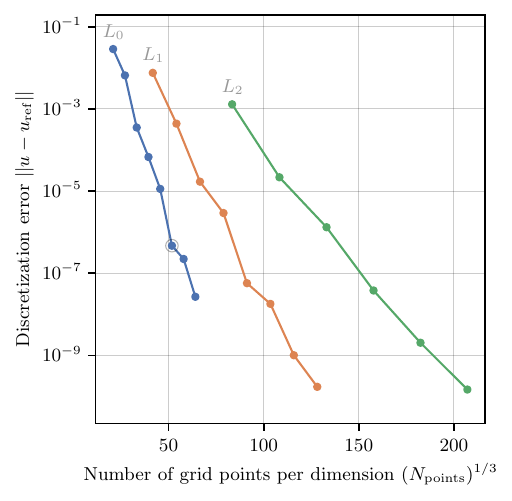}
  \caption{
    \label{fig:bbh_convergence}
    Exponential convergence of the three-dimensional black-hole binary problem
    under uniform $p$~refinement (solid lines) for three uniform $h$~refinement
    levels. The circled configuration is pictured in \cref{fig:bbh_domain}.
    Plotted here is the $L_2$~error~\eqref{eq:errnorm_ref} over all variables
    $\{\psi,\alpha\psi,\beta_\mathrm{excess}^i\}$ and all three interpolation
    points $\bm{x}_m$.}
\end{figure}

Since the binary black hole problem has no analytic solution we assess the
precision of numerical solutions by comparing them to a high-resolution
reference configuration. Specifically, we interpolate all five fields $u_A =
\{\psi, \alpha\psi, \beta_\mathrm{excess}^i\}$ to a set of sample points
$\bm{x}_m$. Then, we compute the discretization error as an $L_2$~norm of the
difference to the high-resolution reference run over all fields and sample
points,
\begin{equation}\label{eq:errnorm_ref}
  \lVert u - u_\mathrm{ref} \rVert \defeq \left(
    \sum_{A,m} \left(u_A(\bm{x}_m) - u_{A,\mathrm{ref}}(\bm{x}_m)\right)^2
  \right)^{1/2}
  \mkern-22mu
  \text{.}
\end{equation}

\begin{table}
\begin{ruledtabular}
\begin{tabular}{lrrr}
   & $\bm{x}_1 = (8.846,0,0)$ & $\bm{x}_2 = (0,0,0)$ & $\bm{x}_3 = (100,0,0)$ \\
  \hline
  $\psi$ & $1.0919141$ & $1.0602545$ & $1.0033643$ \\
  $\alpha\psi$ & $0.7072066$ & $0.9381658$ & $0.9966282$ \\
  $\beta_\mathrm{excess}^x$ & $0.3870172$ & $0$ & $0.0008802$ \\
  $\beta_\mathrm{excess}^y$ & $-0.1273493$ & $0$ & $-0.0003467$ \\
  $\beta_\mathrm{excess}^z$ & $0$ & $0$ & $0$ \\
\end{tabular}
\caption{
  \label{tab:bbh_intrp_values}
  Sample points $\bm{x}_m$ used in \eqref{eq:errnorm_ref} and
  the value of the reference solution at the sample points.}
\end{ruledtabular}
\end{table}

\Cref{fig:bbh_convergence} presents the convergence of the discretization error
under uniform $hp$~refinement for our strong compact DG
scheme~\eqref{eq:dg_op_compact_strong} with $\penaltyparam=1$. Specifically, we
obtain $h$-refinement levels from the domain depicted in \cref{fig:bbh_domain} by
splitting all elements in two along their three logical axes. We obtain $p$-refinement
levels by incrementing the number of grid points by one in all elements and
dimensions. The DG scheme recovers exponential convergence under $p$~refinement,
and suggests the same $\mathcal{O}(h^P)$~convergence under $h$~refinement that we
have found for the single black hole problem in \cref{sec:kerrschild}. We have
chosen $M_n=0.4229$, $\bm{C}_n=(\pm 8, 0, 0)$, $\Omega_0=0.0144$, $w_n=4.8$,
$R=300$, and sample points along the $x$~axis at $x_1=8.846$ (near horizon),
$x_2=0$ (origin) and $x_3=100$ (far field) here. For the high-resolution
reference configuration in \cref{eq:errnorm_ref} we use a run that is $h$~refined
twice, and has one grid point more per element and dimension than the
highest-resolution configuration included in \cref{fig:bbh_convergence}. The
reference values at the interpolation points are listed in
\cref{tab:bbh_intrp_values}. We have verified that these values are consistent
with the same problem solved with the \spec{}~\cite{Pfeiffer2003-mt, spec}
code up to an absolute error of at most $10^{-7}$, which is the precision we
report in \cref{tab:bbh_intrp_values}.  

In forthcoming work we intend to employ the DG scheme that we have presented
here to develop a scalable initial-data solver for binaries involving black
holes and neutron stars in the \spectre{} numerical relativity
code~\cite{ellsolver}.

\section{Conclusion and future work}\label{sec:conclusion}

We have presented a unified discontinuous Galerkin (DG) internal-penalty scheme
that is applicable to a wide range of elliptic equations. Our scheme applies to
linear and nonlinear second-order elliptic PDEs of one or more variables, where
the variables can be scalars, vectors, or tensors of higher rank. It does not
require problem-specific modifications of the DG discretization or of the
numerical fluxes that couple neighboring elements. The scheme supports a wide
range of linear and nonlinear boundary conditions, and applies to equations
formulated on curved manifolds. We demonstrate its versatility by solving a
simple Poisson problem, a linear elasticity problem on a curved mesh with
nonconforming element boundaries, and two nonlinear problems in general
relativity involving black holes. The unified DG scheme is capable of solving
these problems with no structural changes. It recovers optimal
$\mathcal{O}(h^{P+1})$~convergence for the linear test problems and
$\mathcal{O}(h^P)$~convergence for the nonlinear test problems, where $P$ is
the polynomial degree of the elements. The scheme is implemented in the
open-source \spectre{} code~\cite{spectre} and the results presented in this
article are reproducible with the supplemental input-file
configurations.\footnote{\url{https://arxiv.org/src/2108.05826/anc}}

The DG scheme developed here can potentially be improved in multiple ways in
future work. Dealiasing techniques have the potential to increase the accuracy
of the scheme on curved meshes and for equations with background quantities. The
choice of penalty on curved meshes remains a subject of ongoing study.
Furthermore, detailed studies of the symmetry of the DG operator and related
adjustments to the scheme, such as switching to the strong-weak formulation, can
potentially make the DG operator faster to solve.

Since the convergence properties of the DG scheme are sensitive to the specifics
of the computational domain, we have chosen to refine the domains as uniformly
as possible while retaining some important features, such as curved meshes and
nonconforming element boundaries. For practical applications it is typically
more important to obtain steep rather than uniform convergence, in order to
conserve computational resources and thus achieve faster or more precise solves.
Therefore, a focus of future work will be to develop adaptive mesh-refinement
strategies for the elliptic DG scheme that place grid points in regions and
dimensions of the domain that dominate the discretization error.

Once the DG discretization of the elliptic equations is at hand, numerical
techniques for solving the resulting matrix equation become important.
Sophisticated linear and nonlinear iterative algorithms are necessary to solve
high-resolution elliptic problems in parallel on large computing clusters. Many
of the choices we have made in the development of the DG scheme are motivated by
such large-scale applications. For this purpose we are developing a scalable
multigrid-Schwarz preconditioned Newton-Krylov iterative solver with task-based
parallelism that will be presented in \cite{ellsolver}.

\begin{acknowledgments}
The authors thank Trevor Vincent, Hannes R{\"u}ter, Nils Deppe, Will Throwe,
Lawrence Kidder, and Saul Teukolsky for helpful discussions. N.\ F.\ also thanks the
Cornell Center for Astrophysics and Planetary Science and TAPIR at Caltech for
the hospitality and financial support during research stays. Computations were
performed with the \spectre{} code~\cite{spectre} on the Minerva cluster at the
Max Planck Institute for Gravitational Physics. The figures in this article were
produced with \texttt{dgpy}~\cite{dgpy},
\texttt{matplotlib}~\cite{matplotlib1,matplotlib2}, \texttt{TikZ}~\cite{tikz},
and \texttt{ParaView}~\cite{paraview}.
\end{acknowledgments}

\appendix*

\section{Physical systems}\label{app:systems}

\subsection{Puncture equation}

A popular approach to produce initial data for general-relativistic time
evolutions involving black holes is to reduce the Einstein constraint equations
to a single nonlinear elliptic PDE, the \emph{puncture equation} \footnote{See,
e.g., Section~12.2 in \ccite{BaumgarteShapiro} for an introduction to puncture
initial data.}
\begin{equation}\label{eq:punctures}
  -\partial_i\partial_i u = \beta \left(\alpha \left(1 + u\right) + 1\right)^{-7}
  \text{,}
\end{equation}
written here in Cartesian coordinates. The puncture
equation~\eqref{eq:punctures} is solved for the field~$u(\bm{x})$ from which an
admissible spacetime metric can be constructed \cite{BaumgarteShapiro}. The
quantities
\begin{subequations}\label{eq:punctures_background}  
\begin{align}
  \frac{1}{\alpha} &= \sum_I \frac{M_I}{r_I}
  \text{,} \\
  \beta &= \frac{1}{8}\alpha^7\bar{A}_{ij}\bar{A}^{ij}
  \text{,} \\
  \bar{A}^{ij} &= \frac{3}{2}
  \sum_I \frac{1}{r_I^2} \bigg(\begin{aligned}[t]
    &2 P_I^{(i} n_I^{j)} - (\delta^{ij} - n_I^i n_I^j) P_I^k n_I^k \\
    &+ \frac{4}{r_I}n_I^{(i}\epsilon^{j)kl}S^I_k n^I_l \bigg)
  \end{aligned}
\end{align}
\end{subequations}
are background fields that define a configuration of multiple black holes. The
black holes are parametrized by their puncture masses~$M_I$,
positions~$\bm{C}_I$, linear momenta~$\bm{P}_I$ and spins~$\bm{S}_I$.
In \cref{eq:punctures_background}, $r_I=\lVert \bm{x} - \bm{C}_I \rVert$ is the
Euclidean coordinate distance to the $I$th black hole and
$\bm{n}_I=(\bm{x}-\bm{C}_I) / r_I$ is the radial unit normal to the $I$th black
hole.

To formulate the puncture equation~\eqref{eq:punctures} in first-order
flux form~\eqref{eq:fluxform} we can choose the auxiliary
variable~$v_i=\partial_i u$ and the fluxes and sources
\begin{subequations}\label{eq:punctures_fluxform}
\begin{alignat}{3}
  \tensor{\dgF}{_{\!v}^i_j} &= u \, \delta^i_j
  \text{,} &\quad
  \tensor{\dgS}{_v_{\,j}} &= v_j
  \text{,} \\
  \dgFi{u} &= v_i
  \text{,} &\quad
  \dgS_u &= -\beta \left(\alpha \left(1 + u\right) + 1\right)^{-7}
  \text{,}
\end{alignat}
\end{subequations}
along with both~$f_\alpha=0$. Other Poisson-type equations with nonlinear
sources can be formulated analogously.

\subsection{The XCTS equations of general relativity}\label{sec:xcts}

The XCTS equations
\begin{subequations}\label{eq:xcts}
\begin{align}
  \label{eq:xcts_hamiltonian}
  \bar{\nabla}^2 \psi &= \begin{aligned}[t]
    &\frac{1}{8}\psi\bar{R} + \frac{1}{12}\psi^5 K^2 \\
    &-\frac{1}{8}\psi^{-7}\bar{A}_{ij}\bar{A}^{ij} - 2\pi\psi^5\rho
  \end{aligned}
  \\
  \label{eq:xcts_lapse}
  \bar{\nabla}^2\left(\alpha\psi\right) &= \alpha\psi
  \begin{aligned}[t]
    &\bigg(\frac{7}{8}\psi^{-8}\bar{A}_{ij}\bar{A}^{ij} + \frac{5}{12}\psi^4 K^2 + \frac{1}{8}\bar{R} \\
    &+ 2\pi\psi^4\left(\rho + 2S\right)\bigg) - \psi^5\partial_t K + \psi^5\beta^i\bar{\nabla}_i K
  \end{aligned}
  \\
  \label{eq:xcts_momentum}
  \bar{\nabla}_i(\bar{L}\beta)^{ij} &= \begin{aligned}[t]
    &(\bar{L}\beta)^{ij}\bar{\nabla}_i \ln(\bar{\alpha}) + \bar{\alpha}\bar{\nabla}_i\left(\bar{\alpha}^{-1}\bar{u}^{ij}\right) \\
    &+ \frac{4}{3}\bar{\alpha}\psi^6\bar{\nabla}^j K + 16\pi\bar{\alpha}\psi^{10}S^j
  \end{aligned}
\end{align}
\end{subequations}
with $\bar{A}^{ij} = \frac{1}{2\bar{\alpha}}\left((\bar{L}\beta)^{ij} -
\bar{u}^{ij}\right)$ and $\bar{\alpha} = \alpha \psi^{-6}$ are a set of
nonlinear elliptic equations that the spacetime metric of general relativity
must satisfy at all times~\cite{Pfeiffer2004-oo}.\footnote{See, e.g.,
\ccite{BaumgarteShapiro} for an introduction to the XCTS equations, in
particular Box~3.3.} They are solved for the conformal factor~$\psi$, the
product of lapse and conformal factor~$\alpha\psi$, and the shift
vector~$\beta^j$. The remaining quantities in the equations, i.e., the conformal
metric~$\bar{\gamma}_{ij}$, the trace of the extrinsic curvature~$K$, their
respective time derivatives~$\bar{u}_{ij}$ and $\partial_t K$, the energy
density~$\rho$, the stress-energy trace~$S$ and the momentum density~$S^i$, are
freely specifiable fields that define the scenario at hand. Of particular
importance is the conformal metric~$\bar{\gamma}_{ij}$, which defines the
background geometry, the covariant derivative~$\bar{\nabla}$, the Ricci
scalar~$\bar{R}$ and the longitudinal operator
\begin{equation}
  \left(\bar{L}\beta\right)^{ij} = \bar{\nabla}^i\beta^j + \bar{\nabla}^j\beta^i
  - \frac{2}{3}\bar{\gamma}^{ij}\bar{\nabla}_k\beta^k
  \text{.}
\end{equation}
Note that the XCTS equations are essentially two Poisson equations and one
elasticity equation with coupled, nonlinear sources on a curved manifold. In
this analogy, the longitudinal operator plays the role of the elastic
constitutive relation that connects the symmetric \enquote{shift
strain}~$\bar{\nabla}_{(i}\beta_{j)}$ with the
\enquote{stress}~$(\bar{L}\beta)^{ij}$ of which we take the divergence in the
momentum constraint~\eqref{eq:xcts_momentum}. This particular constitutive
relation is equivalent to an isotropic and homogeneous
material~\eqref{eq:constrel_isotropic_homogeneous} with bulk modulus~$K=0$ (not
to be confused with the extrinsic curvature trace~$K$ in this context) and shear
modulus~$\mu=1$.

To formulate the XCTS equations in first-order flux form~\eqref{eq:fluxform} we
choose for auxiliary variables the gradient of the conformal factor,
$v_i=\partial_i\psi$, the gradient of the lapse times the conformal factor,
$w_i=\partial_i\left(\alpha\psi\right)$, and the symmetric shift strain
$B_{ij}=\bar{\nabla}_{(i}\beta_{j)}$. Then, the XCTS equations~\eqref{eq:xcts}
can be formulated with the fluxes and sources
\begin{widetext}
\begin{subequations}\label{eq:xcts_fluxform}
\begin{alignat}{3}
  \tensor{\dgF}{_{\!v}^i_j} &= \delta^i_j \psi
  \text{,} &\quad
  \tensor{\dgS}{_v_{\, j}} &= v_j
  \text{,} \\
  \tensor{\dgF}{_{\!\psi}^i} &= \bar{\gamma}^{ij} v_j
  \text{,} &\quad
  \dgS_\psi &= \begin{aligned}[t]
    &-\bar{\Gamma}^i_{ij} \tensor{\dgF}{_{\!\psi}^j} + \frac{1}{8}\psi\bar{R} + \frac{1}{12}\psi^5 K^2 
    - \frac{1}{8}\psi^{-7}\bar{A}_{ij}\bar{A}^{ij} - 2\pi\psi^5\rho
  \end{aligned}
\intertext{for \cref{eq:xcts_hamiltonian},}
  \tensor{\dgF}{_{\!w}^i_j} &= \delta^i_j \alpha\psi
  \text{,} &\quad
  \tensor{\dgS}{_w_{\, j}} &= w_j
  \text{,} \\
  \tensor{\dgF}{_{\!(\alpha\psi)}^i} &= \bar{\gamma}^{ij} w_j
  \text{,} &\quad
  \dgS_{(\alpha\psi)} &= \begin{aligned}[t]
    &-\bar{\Gamma}^i_{ij} \tensor{\dgF}{_{\!(\alpha\psi)}^j} 
    + \alpha\psi \bigg(\frac{7}{8}\psi^{-8} \bar{A}_{ij}\bar{A}^{ij} + \frac{5}{12} \psi^4 K^2 
    + \frac{1}{8}\bar{R} + 2\pi\psi^4\left(\rho + 2S\right) \bigg) \\
    &- \psi^5\partial_t K + \psi^5\beta^i\bar{\nabla}_i K
  \end{aligned}
\end{alignat}
for \cref{eq:xcts_lapse}, and
\begin{alignat}{3}
  \tensor{\dgF}{_{\!B}^i_{jk}} &= \delta^i_{(j} \bar{\gamma}_{k)l} \beta^l
  \text{,} &\quad 
  \tensor{\dgS}{_B_{\, j}_k} &= B_{jk} + \bar{\Gamma}_{ijk}\beta^i
  \text{,} \\
  \label{eq:xcts_flux_shift}
  \tensor{\dgF}{_{\!\beta}^i^j} &= 2\big(\bar{\gamma}^{ik}\bar{\gamma}^{jl} - \frac{1}{3} \bar{\gamma}^{ij}\bar{\gamma}^{kl}\big) B_{kl}
  \text{,} &\quad
  \tensor{\dgS}{_\beta^i} &= \begin{aligned}[t]
    &-\bar{\Gamma}^j_{jk} \tensor{\dgF}{_{\!\beta}^i^k} - \bar{\Gamma}^i_{jk} \tensor{\dgF}{_{\!\beta}^j^k} 
    + \left(\tensor{\dgF}{_{\!\beta}^i^j} - \bar{u}^{ij}\right) \bar{\gamma}_{jk} \left(\frac{\tensor{\dgF}{_{\!(\alpha\psi)}^k}}{\alpha\psi} - 7 \frac{\tensor{\dgF}{_{\!\psi}^k}}{\psi}\right) \\
    &+ \bar{\nabla}_j\bar{u}^{ij} + \frac{4}{3}\alpha\bar{\nabla}^i K + 16\pi\alpha\psi^4 S^i
  \end{aligned}
\end{alignat}
\end{subequations}
for \cref{eq:xcts_momentum}. All fixed sources~$\dgf_\alpha(\bm{x})$ vanish.
Note that in \cref{eq:xcts_flux_shift},
$\tensor{\dgF}{_{\!\beta}^i^j}=\left(\bar{L}\beta\right)^{ij}$, expressed in the
auxiliary variables.
\end{widetext}

\bibliography{References}

\end{document}